\documentclass{article} 

\usepackage[english]{babel} 
\usepackage{amssymb}
\usepackage{amsmath}
\usepackage{txfonts}
\usepackage{mathdots}
\usepackage[classicReIm]{kpfonts}
\usepackage{graphicx}
\usepackage{xcolor}

\def\RR{\vbox {\hbox to 8.9pt {I\hskip-2.1pt R\hfil}}}
\def\CC{{\rm C\hskip-4.8pt \vrule height 6pt width 12000sp\hskip 5pt}}

\def\rec#1{\frac{1}{#1}}

\begin{document}




\title{DIFFERENTIATION OF THE WRIGHT FUNCTIONS WITH RESPECT TO PARAMETERS AND OTHER RESULTS}
\maketitle


\vskip -0.6truecm
\noindent \textbf{Alexander Apelblat}               

\noindent Department of Chemical Engineering, Ben Gurion University of the Negev, Beer Sheva, 84105, Israel, {\it apelblat@bgu.ac.il}

\noindent \textbf{}

\noindent \textbf{Francesco Mainardi}

\noindent Department of  Physics and Astronomy, University of Bologna,
Via Irnerio 46, 40126 Bologna, Italy,  {\it francesco.mainardi@bo.infn.it}
\date{}
\noindent \textbf{}
\vskip 0.1truecm

\noindent \textbf{Abstract:} 
In this work we discuss the  derivatives of the Wright functions
(of the first and the second kind) 
 with respect to parameters. Differentiation of these functions leads to infinite power series with coefficient being quotients of the digamma (psi) and gamma functions. Only in few cases it is possible to obtain the sums of these series in a closed form. Functional form of the power series resembles those derived for the Mittag-Leffler functions. If the Wright functions are treated as the generalized Bessel functions, differentiation operations can be expressed in terms of the Bessel functions and their derivatives with respect to the order. It is demonstrated that in many cases it is possible to derive the explicit form of the Mittag-Leffler functions by performing simple operations with the Laplace transforms of the Wright functions. 
 The Laplace transform pairs of the both kinds of the Wright   are discussed for particular values of the parameters. Some transform pairs   serve to obtain functional limits by applying the shifted Dirac delta function.
{
 We expect that the present analysis will find several applications in physics and more generally in applied sciences. 
 Indeed these special functions of the Mittg-Leffler and Wright type  have already found application in rheology and in stochastic processes whre fractional calculus is relevant. Then  careful readers can take benefit from our new results
 presented in this paper for novel applications.  
}
\vskip 0.1truecm


\noindent
\vskip 0.1truecm


\vskip 0.1truecm
\noindent{\bf Published in Appl. Sci. (MDPI), Vol. 2022, No. 12,
12825/1-17 (2022).}
\\
{\bf DOI: 10.3390/app122412825}
\newpage 
\noindent

\noindent

\noindent \textbf{}

\noindent \textbf{}

\section{Introduction}

 Partial differential equations of fractional order are successively applied for modeling time and space diffusion, stochastic processes, probability distributions and other diverse natural phenomena. They are extremely important in physical processes that can be described by using the fractional calculus. In the mathematical literature, when solution of these fractional differential equations is desired, we frequently encounter introduced and named after him, the Wright functions. At beginning, at 1933 [1] and at 1940 [2], these functions were considered as a some kind of generalization of the Bessel functions, but today they play an significant independent role in the theory of special functions. There are many investigations devoted to analytical properties and applications of the Wright functions, but here are mentioned only two survey papers where essential material on the subject is included [3,4].
 These functions turn out to be  particular cases of higher transcendental functions as 
recently shown in interesting surveys by Kiryakova [5],  Srivastava [6] .
 
 In this paper we discuss  three quite different subjects which are associated with the Wright functions. In the first part, the Wright function $W_{\alpha, \beta}(t)$ where $t$ is the argument and $\alpha$ and $\beta$ are the parameters, is differentiated with respect to parameters and derived expressions are compared with similar formulas for the Mittag-Leffler functions. In a continuous effort, after investigating differentiation of the Bessel functions and the Mittag-Leffler functions with respect to their parameters [7-9], this mathematical operation is extended here to the Wright functions. Special attention is devoted to the cases when the Wright functions can be reduced to the Bessel functions and expressed in a closed form. The auxiliary functions, $F_{\alpha}(t)$ and $M_{\alpha}(t)$, which were introduced for the first time in 1990's  by Mainardi, see  [4], and now are called the Mainardi functions, are as well discussed in this section.
 
  Functional behaviour of derivatives with respect to the order is also presented in graphical form. The presented plots were prepared by evaluation of sums of infinite series by using MATHEMATICA program.

 The second part of this paper is dedicated to the Laplace transform pairs of the Wright functions. It is demonstrated how the Laplace transforms of the Wright functions are useful for obtaining explicit expressions for the Mittag-Leffler functions.
 Finally, we  discuss the functional limits which are associated with the Wright and the Mittag-Leffler functions. These limits can be derived by applying the delta sequence in the form of the shifted Dirac function. This delta sequence is directly related to the order of Bessel function and was introduced by Lamborn in 1969, see [8-11] and [18-20].
\newpage

 Throughout this paper all mathematical operations or manipulations with functions, series, integrals, integral representations and transforms are formal and it is assumed that arguments and parameters are real  numbers. There will be no proofs of validity of derived results, though they are presumed  to be correct considering that in a part they were previously obtained independently by other methods.

\noindent

\section{The Wright functions of the first and \\ of the second kind}

The Wright functions $W_{\alpha, \beta}(z)$ are defined by the series representation as a function of complex argument $z$ and parameters $\alpha$ and $\beta$.
\begin{equation} \label{GrindEQ__2_1_}
W_{\alpha ,\beta } (z)=\sum _{k{\kern 1pt} ={\kern 1pt} 0}^{\infty }
\frac{z^{k} }{k!\, \, \Gamma (\alpha {\kern 1pt} {\kern 1pt} k+\beta )}.
\end{equation}
They are entire functions of $z \in \CC $ for $\alpha > -1$ and for any complex 
$\beta$  (here always $\beta \ge 0$).  
According to Mainardi  see i.e. the appendix F of [4],
we  distinguish the Wright function of the first kind for $\alpha \ge 0$, and of the second kind for $- 1 < \alpha < 0$.
This distinctions is justified for the difference in the aymptotics representations in the complex domain and in the Laplace transforms for the  real positive argument.
For our purposes we recall their Laplace transforms for positive argument $t$,
We have, by  using the symbol $\div$ to denote the  juxtaposition of a function $f(t)$ with its Laplace transform $\widetilde f(s)$, 
\\
for the first kind, when $\alpha \ge 0$
\begin{equation}
  W_{\alpha ,\beta } (\pm t) \,\div \,
    \rec{s}\, E_{\lambda ,\mu }\left(\pm \rec{s}\right) \,;
\end{equation}
\\
for the second kind,   when $-1<\alpha<0$
and putting for convenience  $\nu = -\alpha$
so $0< \nu<1$
\begin{equation}
W_{-\nu ,\beta } (-t) \,\div \,
    E_{\nu ,\beta+\nu }\left(-s \right) \,.
\end{equation}
Above we have  introduced   
the  Mittag-Leffler function
in two parameters $\alpha>0$, $\beta \in \CC$ 
  defined as its convergent series for  all  $z \in \CC$
  \begin{equation}
E_{\alpha, \beta}(z) := \sum_{k=0}^{\infty}
{\frac{z^k}{\Gamma(\alpha k + \beta)}}.
\label{eq:mittag-leffler}
\end{equation}
For more details on the special functions of the Mittag-leffler type we refer the interested readers to the treatise by Gorenflo et al [7] 
where in the recent 2-nd edition also the Wright functions are treated 
in some detail. 

\newpage

\section{Differentiation of the Wright functions \\ of the first kind with respect to parameters}

 We first compare  compare the Wright functions of the first kind with  the two-parameter  Mittag-Leffler functions 
 for $\alpha>0$ and $\beta \ge0$
 from which they differ only by the absence of factorials

Direct differentiation of series with respect to the parameter $\alpha$ gives
\begin{equation} \label{GrindEQ__2_3_}
\begin{array}{l} {\frac{\partial {\kern 1pt} W_{\alpha ,\beta } (t)}{\partial \alpha } =-\sum _{k{\kern 1pt} ={\kern 1pt} \, 1}^{\infty }\left(\frac{{\kern 1pt} \psi (\alpha {\kern 1pt} k+\beta )}{k{\kern 1pt} !\, \, \Gamma (\alpha {\kern 1pt} k+\beta )} \right)\, k\, t^{k}  =} \\ {-\sum _{k{\kern 1pt} ={\kern 1pt} \, 1}^{\infty }\left(\frac{{\kern 1pt} \psi (\alpha {\kern 1pt} k+\beta )}{(k-1){\kern 1pt} !\, \, \Gamma (\alpha {\kern 1pt} k+\beta )} \right)\, \, t^{k}  } \\ {\frac{\partial {\kern 1pt} E_{\alpha ,\beta } (t)}{\partial \alpha } =-\sum _{k{\kern 1pt} ={\kern 1pt} 1}^{\infty }\left(\frac{{\kern 1pt} \psi (\alpha {\kern 1pt} k+\beta )}{\, \Gamma (\alpha {\kern 1pt} k+\beta )} \right)\, k\, t^{k}  } \end{array}
\end{equation}
and with respect to the parameter $\beta$
\noindent
\begin{equation}\label{GrindEQ__2_4_}
\begin{array}{l} {\frac{\partial {\kern 1pt} W_{\alpha ,\beta } (t)}{\partial \beta } =-\sum _{k{\kern 1pt} ={\kern 1pt} \, 0}^{\infty }\left(\frac{{\kern 1pt} \psi (\alpha {\kern 1pt} k+\beta )}{k{\kern 1pt} !\, \, \Gamma (\alpha {\kern 1pt} k+\beta )} \right)\, \, t^{k}  } \\ {\frac{\partial {\kern 1pt} E_{\alpha ,\beta } (t)}{\partial \beta } =-\sum _{k{\kern 1pt} =\, {\kern 1pt} 0}^{\infty }\left(\frac{{\kern 1pt} \psi (\alpha {\kern 1pt} k+\beta )}{\Gamma (\alpha {\kern 1pt} k+\beta )} \right)\, \, t^{k}  } \end{array}
\end{equation}
The second derivatives are
\begin{equation} \label{GrindEQ__2_5_}
\begin{array}{l} {\frac{\partial ^{2} W_{\alpha ,\beta } (t)}{\partial \alpha ^{2} } =\sum _{k{\kern 1pt} =\, {\kern 1pt} 1}^{\infty }\left\{\frac{[\psi (\alpha {\kern 1pt} k+1)]^{2} -\psi ^{(1)} {\kern 1pt} (\alpha {\kern 1pt} k+\beta )}{k{\kern 1pt} {\kern 1pt} !\, \, \Gamma (\alpha {\kern 1pt} k+\beta )} \right\} \, k{\kern 1pt} ^{2} t^{k} } \\ {\frac{\partial ^{2} E_{\alpha ,\beta } (t)}{\partial \alpha ^{2} } =\sum _{k{\kern 1pt} =\, {\kern 1pt} 1}^{\infty }\left\{\frac{[\psi (\alpha {\kern 1pt} k+1)]^{2} -\psi ^{(1)} {\kern 1pt} (\alpha {\kern 1pt} k+\beta )}{\Gamma (\alpha {\kern 1pt} k+\beta )} \right\} \, k{\kern 1pt} ^{2} t^{k} } \end{array}
\end{equation}
and
\begin{equation} \label{GrindEQ__2_6_}
\begin{array}{l} {\frac{\partial ^{2} W_{\alpha ,\beta } (t)}{\partial {\kern 1pt} \beta ^{2} } =\sum _{k{\kern 1pt} ={\kern 1pt} \, 0}^{\infty }\left\{\frac{[\psi (\alpha {\kern 1pt} k+\beta )]^{2} -\psi ^{(1)} {\kern 1pt} (\alpha {\kern 1pt} k+\beta )}{k!\, \, \Gamma (\alpha {\kern 1pt} k+\beta )} \right\} \, \, t^{k} } \\ {\frac{\partial ^{2} E_{\alpha ,\beta } (t)}{\, \partial {\kern 1pt} \beta ^{2} } =\sum _{k{\kern 1pt} ={\kern 1pt} \, 0}^{\infty }\left\{\frac{{\kern 1pt} \psi (\alpha {\kern 1pt} k+\beta )]^{2} -\psi ^{(1)} (\alpha {\kern 1pt} k+\beta )}{\Gamma (\alpha {\kern 1pt} k+\beta )} \right\} \, {\kern 1pt} t^{k} } \end{array}
\end{equation}

We note that for the Mittag-Leffler and the Wright  functions we have the same functional expressions, but in the case of the Wright functions  factorials
appear.
 Contrary to the Mittag-Leffler functions [9-10], summation of these series by using MATHEMATICA gives only few results in a closed form in terms of the generalized hypergeometric functions
\begin{equation} \label{GrindEQ__2_7_}
\begin{array}{l} {\frac{\partial {\kern 1pt} W_{\alpha ,\beta } (t)}{\partial \alpha } \left|{}_{\alpha \, =\, 1,\, \beta \, =\, 0} \right. =-\sum _{k{\kern 1pt} =\, {\kern 1pt} 1}^{\infty }\left(\frac{{\kern 1pt} \psi ({\kern 1pt} k)}{[(k-1)!\, ]^{2} \, \, } \right)\, \, t^{k}  =t_{0} F_{1} (\, \, ;1;t)} \\ {\frac{\partial {\kern 1pt} W_{\alpha ,\beta } (t)}{\partial \alpha } \left|{}_{\alpha \, =\, 1,\, \beta \, =\, 1} \right. =-\sum _{k{\kern 1pt} =\, {\kern 1pt} 1}^{\infty }\left(\frac{{\kern 1pt} \psi ({\kern 1pt} k+1)}{(k-1)!\, k{\kern 1pt} !\, \, } \right)\, \, t^{k}  =t_{0} F_{1} (\, \, ;2;t)} \end{array}
\end{equation}
and
\begin{equation} \label{GrindEQ__2_8_}
\begin{array}{l} {\frac{\partial {\kern 1pt} W_{\alpha ,\beta } (t)}{\partial \beta } \left|{}_{\alpha \, =\, 1,\, \beta \, =\, 0} \right. =-\sum _{k{\kern 1pt} =\, {\kern 1pt} 0}^{\infty }\left(\frac{{\kern 1pt} \psi ({\kern 1pt} k)}{(k-1)!\, k!\, \, \, } \right)\, \, t^{k}  =} \\ {\frac{1}{2} \, [t{\kern 1pt} _{0} F_{1} (\, \, ;1;t)-t{\kern 1pt} _{0} F_{1} (\, \, ;2;t)\, \ln t]+\sqrt{t} K_{1} (2\sqrt{t} )} \\ {\frac{\partial {\kern 1pt} W_{\alpha ,\beta } (t)}{\partial \beta } \left|{}_{\alpha \, =\, 1,\, \beta \, =\, 1} \right. =-\sum _{k{\kern 1pt} =\, {\kern 1pt} 0}^{\infty }\left(\frac{{\kern 1pt} \psi ({\kern 1pt} k+1)}{(\, k{\kern 1pt} !\, )^{2} \, } \right)\, \, t^{k}  =\, _{0} F_{1} (\, \, ;1;t)} \end{array}
\end{equation}
 In the last case, $\alpha = \beta = 1$, in the Brychkov compilation of infinite series [12], the sum is expressed in terms of the modified Bessel functions
\begin{equation} \label{GrindEQ__2_9_}
\begin{array}{l} {\frac{\partial {\kern 1pt} W_{\alpha ,\beta } (t)}{\partial \beta } \left|{}_{\alpha \, =\, \beta \, =\, 1} \right. =-\sum _{k{\kern 1pt} ={\kern 1pt} 0}^{\infty }\left(\frac{{\kern 1pt} \psi ({\kern 1pt} k+1)}{(k{\kern 1pt} !\, )^{2} } \right)\, \, t^{k}  =} \\ {-\frac{1}{2} \, \ln t\, I_{0} (2\sqrt{t} )-K_{0} (2\sqrt{t} )} \end{array}
\end{equation}
 Using MATHEMATICA program, values of derivatives with respect to parameters $\alpha$ and $\beta$ of the Wright functions of the first kind were calculated for the argument $0.25 \le t \le 4.0$ and for parameters $0 \le \alpha \le 5.0$ and $0 \le \beta \le 2.0$. 
 
 In Figure 1 is illustrated behaviour of derivatives with respect to parameter $\alpha$ at different values if the argument $t$. 

 As can be observed, in the $0 < \alpha < 1$ region, exists minimum, with increasing $\alpha$ all curves tend to zero value. The absolute value of the minimum increases with increasing argument.
 
 Derivatives with respect to the parameters $\alpha$ and $\beta$ when the argument $t$ is constant, are presented in Figure 2. The functional form of curves with the change of $\beta$ values is similar to that observed previously in Figure 1.

In order to compare the  behaviour of derivatives with respect to $\alpha $ with those with respect to $\beta$, the same conditions were imposed on $t$ and $\beta$ in Figures 3 and 4 as are in Figures 1 and 2.

 As  can be observed, the similarity of corresponding curves is evident, the only difference is that the absolute values of the minima are lower for derivatives with respect to parameter $\beta$ than to for $\alpha$.

\begin{figure}[h!]
  \includegraphics[width=12cm]{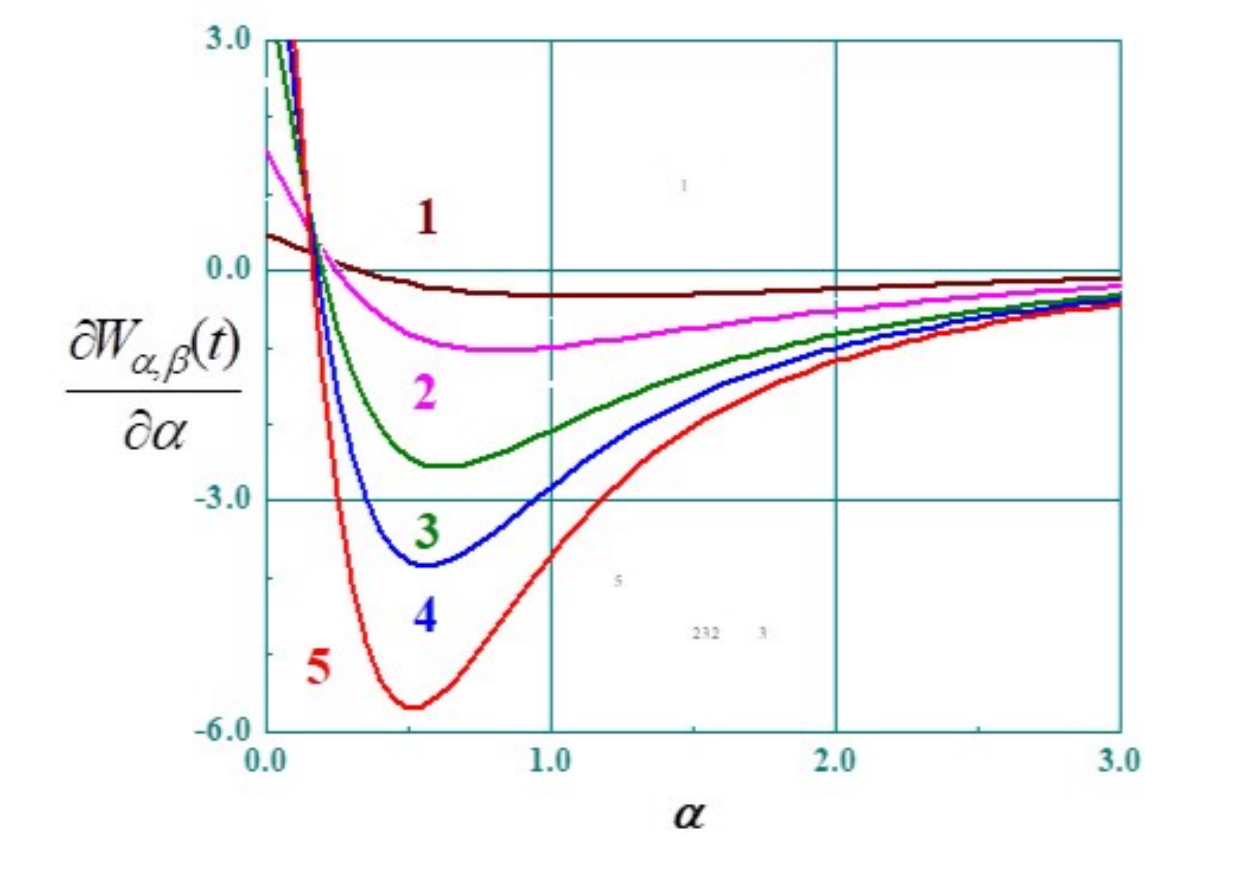}
\vskip -0.2truecm
\caption{Derivatives of the Wright functions of the first kind with respect to parameter $\alpha$ as a function of $\alpha$ for $\beta = 1$ and 
{\bf 1}: $t = 0.5$; {\bf 2}: $t = 1.0$; {\bf 3}: $t = 1.5$;
{\bf 4}: $t = 1.75$;  {\bf 5}: $t = 2.0$.}
\label{Figure: Numerical difference}       
\end{figure}
\begin{figure}[h!]
  \includegraphics[width=12cm]{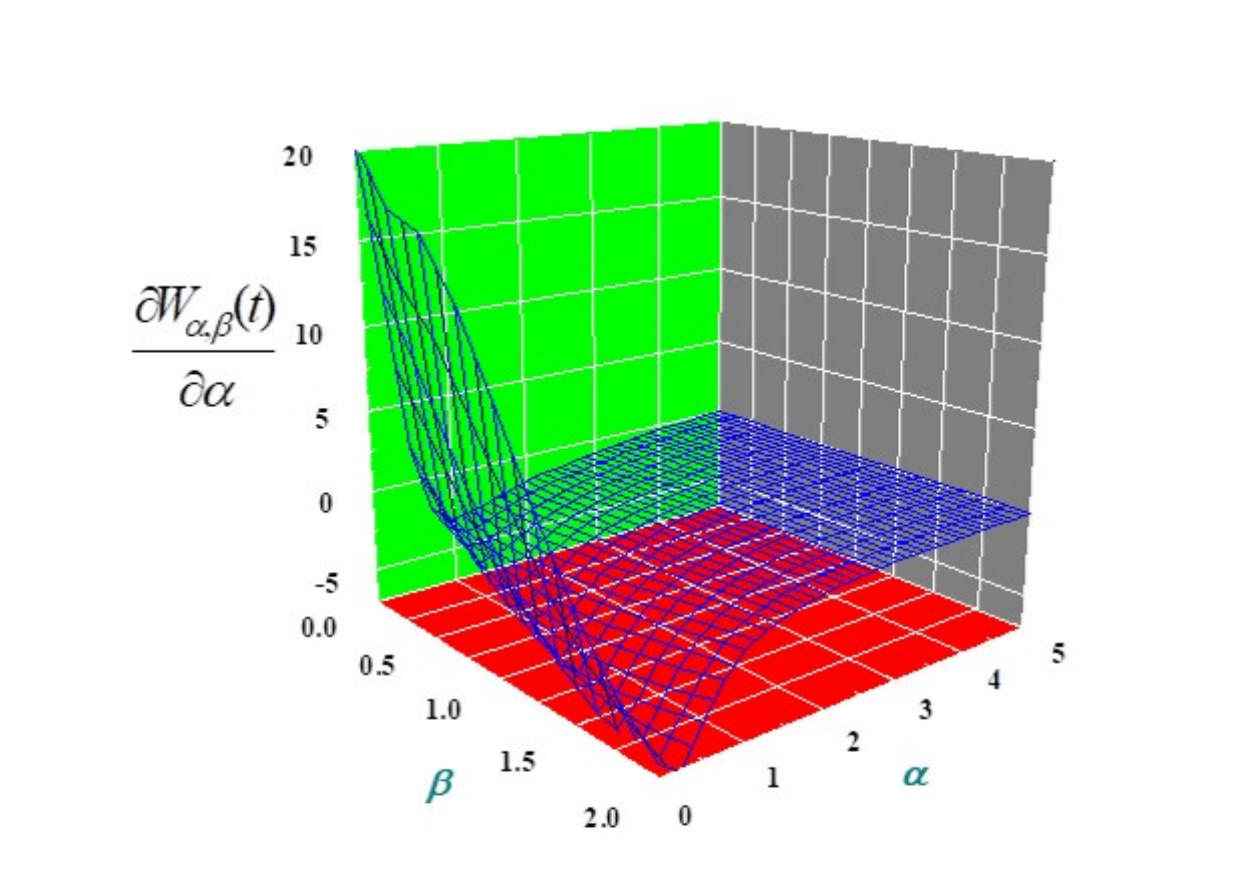}
\vskip -0.2truecm
\caption{Derivatives of the Wright functions of the first kind with respect to parameter $\alpha$ as a function of $\alpha$ and $\beta = 1$ for $t = 2.0$.}
\label{Figure: Numerical difference}       
\end{figure}




\begin{figure}[h!]
  \includegraphics[width=12cm]{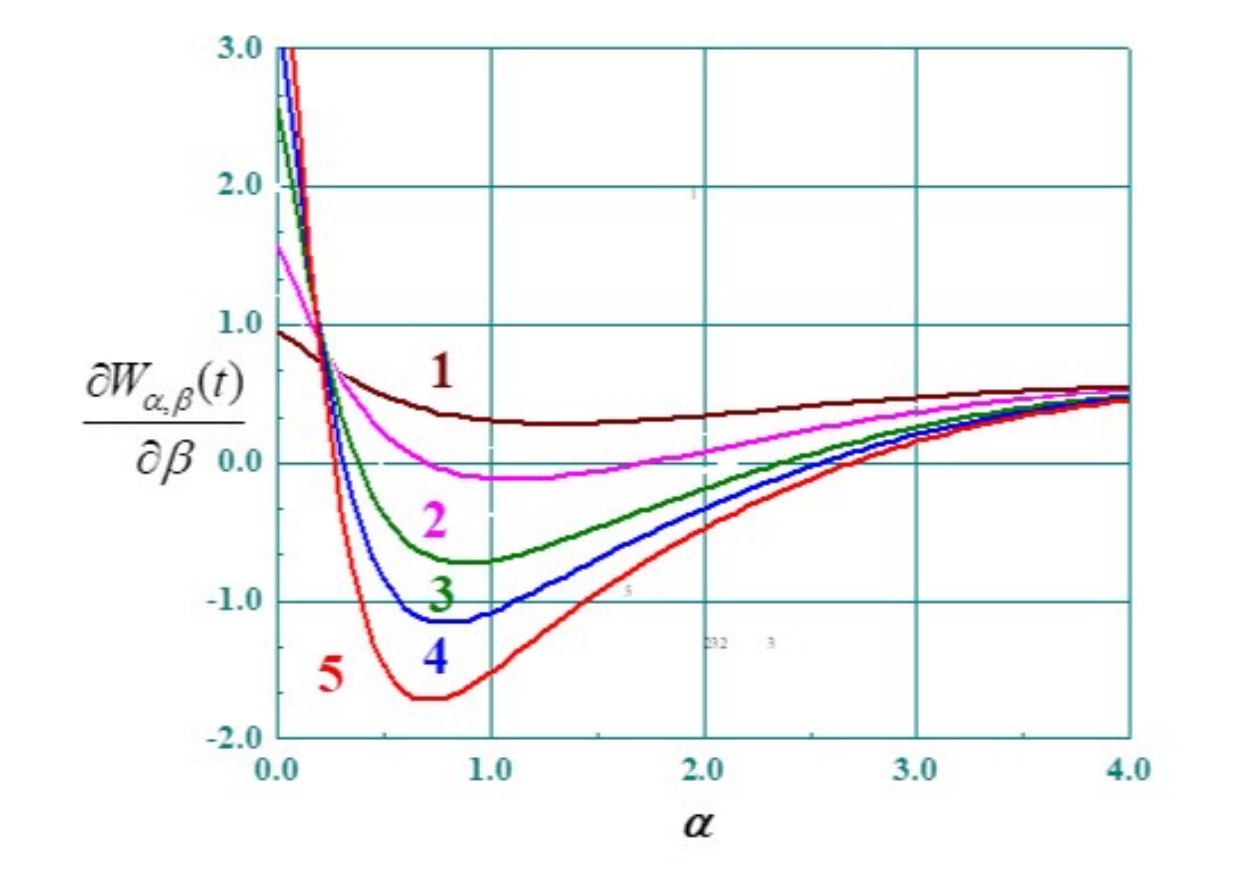}
\vskip -0.2truecm
\caption{Derivatives of the Wright functions with respect to parameter $\beta$ as a function of $\alpha$ for $\beta = 1$ and 
{\bf 1}: $t = 0.5$; {\bf 2} $t = 1.0$;
{ \bf 3} $t = 1.5$; {\bf 4}:$ t = 1.75$; {\bf 5} $t = 2.0$.}
\label{Figure: Numerical difference}       
\end{figure}
\begin{figure}[h!]
  \includegraphics[width=12cm]{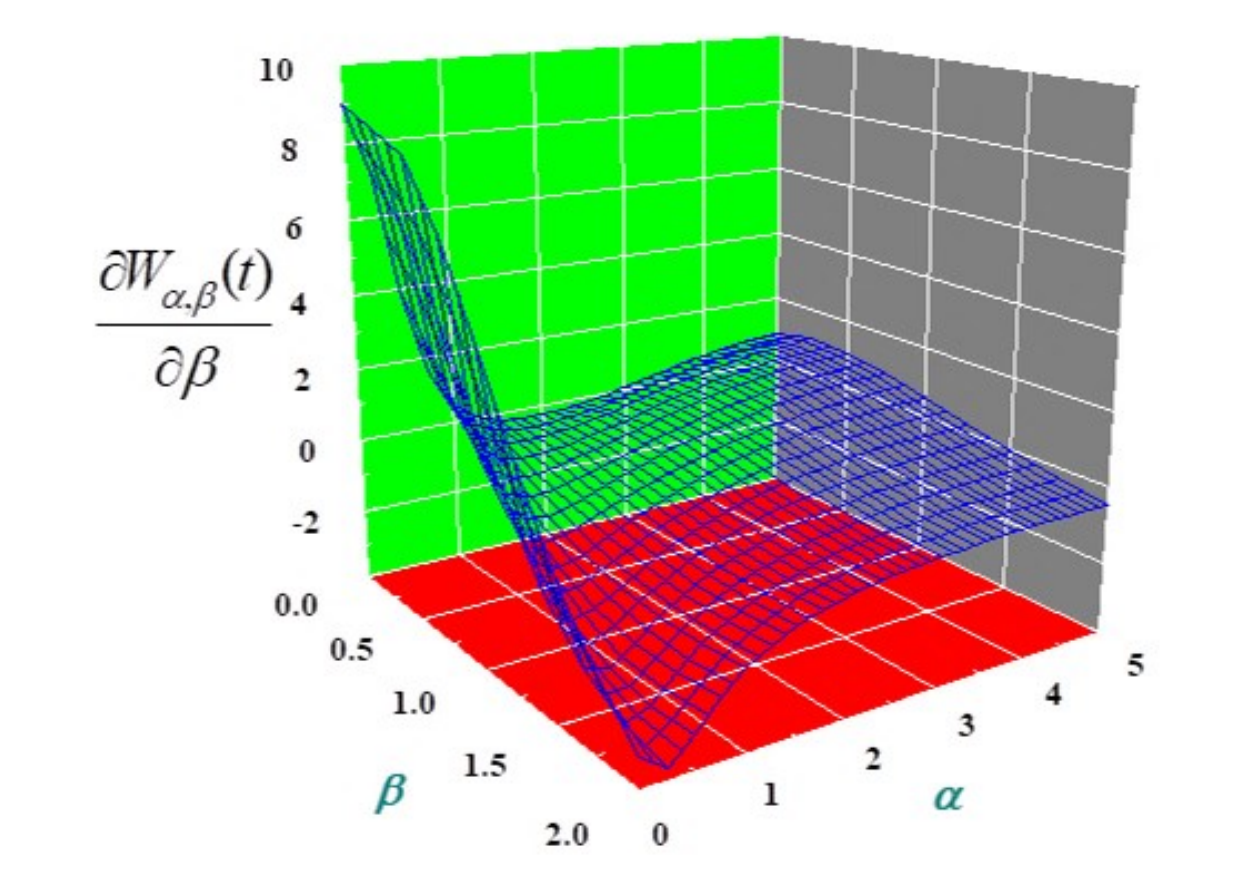}
\vskip -0.2truecm
\caption{Derivatives of the Wright functions with respect to parameter $\beta$ as a function of $\alpha$ and $\beta = 1$ for $t = 2.0$.}
\label{Figure: Numerical difference}       
\end{figure}
%


\noindent
\section{Differentiation of the Wright functions \\ of the  second kind 
with respect to parameters}
 
 We then consider among the Wright functions of the second kind the  functions, $F_{\alpha}(t)$ and $M_{\alpha}(t)$,
 introduced by Mainardi, 
\begin{equation} \label{GrindEQ__2_10_}
\begin{array}{l} {F_{\alpha } (t)=W_{-\, \alpha ,0} (t)\quad ;\quad 0<\alpha <1} \\ {M_{\alpha } (t)=W_{-\, \alpha ,1\, -} \, _{\alpha } (t)\quad ;\quad 0<\alpha <1} \\ {F_{\alpha } (t)=\alpha \, t\, M_{\alpha } (t)} \end{array}
\end{equation}
Their series expansions explicitly read
\begin{equation} \label{GrindEQ__2_11_}
F_{\alpha } (t)=
\sum _{k = 1}^{\infty }\frac{(-\, t)^{k} }{k!\, \Gamma (-\, \alpha \, k)} 
 =-\frac{1}{\pi } \sum _{k = 1}^{\infty }\frac{(-\ t)^{k}}{k!}  \, \Gamma (\alpha \, k+1)\, \sin (\pi \, \alpha \, k)
\end{equation}
and
\begin{equation} \label{GrindEQ__2_12_}
M_{\alpha } (t)=\sum _{k\, =\, 0}^{\infty }\frac{(-\, t)^{k\, } }{k!\, \Gamma \, \left(-\, \alpha (\, k+1)+1\right)}  =\frac{1}{\pi } \sum _{k\, =\, 1}^{\infty }\frac{(-\, t)^{k\, -\, 1} }{(k-1)!\, }  \, \Gamma (\alpha \, k)\, \sin (\pi \, \alpha \, k)
\end{equation}
Direct differentiation of \eqref{GrindEQ__2_11_} and \eqref{GrindEQ__2_12_} gives
\begin{equation} \label{GrindEQ__2_13_}
\begin{array}{l} {\frac{\partial {\kern 1pt} F_{\alpha } (t)}{\partial \alpha } =\frac{1}{\pi } \sum _{k\, =\, 1}^{\infty }\frac{k\, (-\, t)^{k\, -\, 1} }{k!\, }  \, \Gamma (\alpha \, k+1)\, \left[\psi (\alpha \, k+1)\sin (\pi \, \alpha \, k)+\pi \cos (\pi \, \alpha \, k)\right]} \\ {\frac{\partial {\kern 1pt} M_{\alpha } (t)}{\partial \alpha } =\frac{1}{\pi } \sum _{k\, =\, 1}^{\infty }\frac{k\, (-\, t)^{k\, -\, 1} }{(k-1)!\, }  \, \Gamma (\alpha \, k)\, \left[\psi (\alpha \, k)\sin (\pi \, \alpha \, k)+\pi \cos (\pi \, \alpha \, k)\right]} \end{array}
\end{equation}
Using the last equation in \eqref{GrindEQ__2_10_} we have

\noindent
\begin{equation} \label{GrindEQ__2_14_}
\frac{\partial {\kern 1pt} F_{\alpha } (t)}{\partial \alpha } =t\, M_{\alpha } (t)+\alpha \, t\, \frac{\partial {\kern 1pt} M_{\alpha } (t)}{\partial \alpha }
\end{equation}
The second derivatives of these functions are
\begin{equation} \label{GrindEQ__2_15_}
\begin{array}{l} {\frac{\partial {\kern 1pt} ^{2} F_{\alpha } (t)}{\partial \alpha ^{2} } =\frac{1}{\pi } \sum _{k\, =\, 1}^{\infty }\frac{k\, (-\, t)^{k\, -\, 1} }{(k-1)!\, }  \, \Gamma (\alpha \, k+1)\left\{[\psi '(\alpha \, k+1)+\right. \psi ^{2} (\alpha \, k+1)]\sin (\pi \, \alpha \, k)+} \\ {2{\kern 1pt} {\kern 1pt} \pi \cos (\pi \, \alpha \, k)\, \psi (\alpha \, k+1)-\left. \pi ^{2} \sin (\pi \, \alpha \, k)\right\}} \end{array}
\end{equation}
and
\begin{equation} \label{GrindEQ__2_16_}
\begin{array}{l} {\frac{\partial {\kern 1pt} ^{2} M_{\alpha } (t)}{\partial \alpha ^{2} } =\frac{1}{\pi } \sum _{k\, =\, 1}^{\infty }\frac{k^{2} \, (-\, t)^{k\, -\, 1} }{(k-1)!\, }  \, \Gamma (\alpha \, k)\left\{[\psi '(\alpha \, k)+\right. \psi ^{2} (\alpha \, k)]\sin (\pi \, \alpha \, k)+} \\ {2{\kern 1pt} {\kern 1pt} \pi \cos (\pi \, \alpha \, k)\, \psi (\alpha \, k)-\left. \pi ^{2} \sin (\pi \, \alpha \, k)\right\}} \end{array}
\end{equation}
They are interrelated by
\begin{equation} \label{GrindEQ__2_17_}
\frac{\partial ^{2} {\kern 1pt} F_{\alpha } (t)}{\partial \alpha ^{2} } =2\, t\, \frac{\partial {\kern 1pt} M_{\alpha } (t)}{\partial \alpha } +\alpha {\kern 1pt} t\frac{\partial ^{2} {\kern 1pt} M_{\alpha } (t)}{\partial \alpha ^{2} }
\end{equation}



\section{The Laplace transforms and the Wright functions of the first kind}

 The Laplace transforms of the Wright functions are expressed in terms of two-parameter the Mittag-Leffler functions [3,4]
\begin{equation} \label{GrindEQ__3_1_}
L\left\{W_{\alpha ,\beta } (\pm \lambda {\kern 1pt} t)\right\}=\frac{1}{s} \, E_{\alpha ,\beta } \left(\pm \frac{\lambda }{s} \right)\quad ;\quad \alpha >0\quad ;\quad \lambda >0
\end{equation}
Applying operational rules of the Laplace transformation we have, see  [13], [15-16],
\begin{equation} \label{GrindEQ__3_2_}
L\left\{e^{\pm \, \rho } W_{\alpha ,\beta } (\lambda {\kern 1pt} t)\right\}=\frac{1}{s\mp \rho } \, E_{\alpha ,\beta } \left(\frac{\lambda }{s\mp \rho } \right)\quad ;\quad \lambda ,\rho >0
\end{equation}
and this permits to obtain
\begin{equation} \label{GrindEQ__3_3_}
\begin{array}{l} {L\left\{\sinh (\rho {\kern 1pt} t){\kern 1pt} \, W_{\alpha ,\beta } (\lambda {\kern 1pt} t)\right\}=\frac{1}{2} \left\{\frac{1}{s-\rho } \, E_{\alpha ,\beta } \left(\frac{\lambda }{s-\rho } \right)-\frac{1}{s+\rho } \, E_{\alpha ,\beta } \left(\frac{\lambda }{s+\rho } \right)\right\}\quad } \\ {L\left\{\cosh (\rho {\kern 1pt} t){\kern 1pt} \, W_{\alpha ,\beta } (\lambda {\kern 1pt} t)\right\}=\frac{1}{2} \left\{\frac{1}{s-\rho } \, E_{\alpha ,\beta } \left(\frac{\lambda }{s-\rho } \right)+\frac{1}{s+\rho } \, E_{\alpha ,\beta } \left(\frac{\lambda }{s+\rho } \right)\right\}\quad } \end{array}
\end{equation}
Multiplication \eqref{GrindEQ__3_1_} by \textit{t} gives
\begin{equation} \label{GrindEQ__3_4_}
L\left\{t\, W_{\alpha ,\beta } (\lambda {\kern 1pt} t)\right\}=-\frac{d}{ds} \left\{\frac{1}{s} \, E_{\alpha ,\beta } \left(\frac{\lambda }{s} \right)\right\}=-\left\{-\frac{1}{s^{2} } E_{\alpha ,\beta } \left(\frac{\lambda }{s} \right)+\frac{1}{s} \frac{d}{ds} E_{\alpha ,\beta } \left(\frac{\lambda }{s} \right)\right\}
\end{equation}
Derivative of the Mittag-Leffler function is
\begin{equation} \label{GrindEQ__3_5_}
\frac{d}{ds} E_{\alpha ,\beta } \left(\frac{\lambda }{s} \right)=-\frac{\lambda }{s^{2} } \left\{\frac{E_{\alpha ,\beta \, -\, 1} \left(\frac{\lambda }{s} \right)-(\beta -1)E_{\alpha ,\beta } \left(\frac{\lambda }{s} \right)}{\alpha \left(\frac{\lambda }{s} \right)} \right\}
\end{equation}
and finally we have
\begin{equation} \label{GrindEQ__3_6_}
L\left\{t\, W_{\alpha ,\beta } (\lambda {\kern 1pt} t)\right\}=\frac{1}{s^{2} } \left\{\frac{(\alpha {\kern 1pt} \lambda -\beta +1)E_{\alpha ,\beta } \left(\frac{\lambda }{s} \right)+E_{\alpha ,\beta \, -\, 1} \left(\frac{\lambda }{s} \right)}{\alpha {\kern 1pt} \lambda } \right\}
\end{equation}

 In case that the Wright functions are expressed as the Bessel functions (see (2.18)), the Laplace transforms are known for $\beta = 0,1,2$ [12]
\begin{equation} \label{GrindEQ__3_7_}
\begin{array}{l} {\int _{0}^{\infty }e^{-\, s{\kern 1pt} t}  W_{1,1} (-\frac{\lambda ^{2} t^{2} }{4} )\, dt=\int _{0}^{\infty }e^{-\, s{\kern 1pt} t}  \, J_{0} (\lambda {\kern 1pt} t)\, dt=\frac{1}{\sqrt{s^{2} +\lambda ^{2} } } } \\ {\int _{0}^{\infty }e^{-\, s{\kern 1pt} t}  W_{1,2} (-\frac{\lambda ^{2} t^{2} }{4} )\, dt=\frac{2}{\lambda {\kern 1pt} } \int _{0}^{\infty }e^{-\, s{\kern 1pt} t}  \frac{J_{1} (\lambda {\kern 1pt} t)}{t} \, dt=\frac{2{\kern 1pt} {\kern 1pt} }{{\kern 1pt} \, {\kern 1pt} [s+\sqrt{s^{2} +\lambda ^{2} } ]} } \\ {\int _{0}^{\infty }e^{-\, s{\kern 1pt} t}  W_{1,3} (-\frac{\lambda ^{2} t^{2} }{4} )\, dt=\frac{4}{\lambda ^{2} {\kern 1pt} } \int _{0}^{\infty }e^{-\, s{\kern 1pt} t}  \, \frac{J_{2} (\lambda {\kern 1pt} t)}{t^{2} } \, dt=} \\ {\frac{1}{\lambda } \left\{\frac{{\kern 1pt} \lambda }{{\kern 1pt} {\kern 1pt} [s+\sqrt{s^{2} +\lambda ^{2} } ]} +\frac{1}{3} \left[\frac{{\kern 1pt} \lambda }{{\kern 1pt} {\kern 1pt} [s+\sqrt{s^{2} +\lambda ^{2} } ]} \right]^{3} \right\}} \end{array}
\end{equation}
From \eqref{GrindEQ__2_18_} it follows that
\begin{equation}
\int _{0}^{\infty }e^{-\, s{\kern 1pt} t}  W_{1,\beta \, +\, 1} (-\lambda {\kern 1pt} t)\, dt=\frac{1}{\lambda ^{\beta /2} } \, \int _{0}^{\infty }e^{-\, s{\kern 1pt} t}  t^{-\, \beta /2} J_{\beta } (2\sqrt{\lambda {\kern 1pt} t} )\, dt
\label{GrindEQ__3_8_}
\end{equation}
and this integral equality is useful to derive explicit forms of the Mittag-Leffler functions. Starting with $\beta = 0$ we have [14]
\begin{equation} \label{GrindEQ__3_9_}
\int _{0}^{\infty }e^{-\, s{\kern 1pt} t}  W_{1,1} (-\lambda {\kern 1pt} t)\, dt=\, \int _{0}^{\infty }e^{-\, s{\kern 1pt} t}  J_{0} (2\sqrt{\lambda {\kern 1pt} t} )\, dt=\frac{1}{s} \, e^{-\, \lambda /s} \quad ;\quad Res>0
\end{equation}
but from \eqref{GrindEQ__3_1_}
\begin{equation} \label{GrindEQ__3_10_}
L\left\{W_{1,1} (-\lambda {\kern 1pt} t)\right\}=\frac{1}{s} \, E_{1,1} (-\frac{\lambda }{s} )
\end{equation}
and therefore by comparing the expected result is reached
\begin{equation} \label{GrindEQ__3_11_}
\begin{array}{l} {E_{1,1} (-\frac{\lambda }{s} )=e^{-\, \lambda /s} } \\ {\tau =\frac{\lambda }{s} } \\ {E_{1,1} (-\, \tau )=e^{-\, \tau } } \end{array}
\end{equation}
Introducing $\beta = 1$ into \eqref{GrindEQ__3_8_} we have [14]
\begin{equation} \label{GrindEQ__3_12_}
\begin{array}{l} {\int _{0}^{\infty }e^{-\, s{\kern 1pt} t}  W_{1,2} (-\lambda {\kern 1pt} t)\, dt=\, \frac{1}{\sqrt{\lambda } } \int _{0}^{\infty }e^{-\, s{\kern 1pt} t}  t^{-\, 1/2} J_{1} (2\sqrt{\lambda {\kern 1pt} t} )\, dt=} \\ {\sqrt{\frac{\pi }{\lambda {\kern 1pt} s} } \, e^{-\, \lambda /2{\kern 1pt} s} I_{1/2} (\frac{\lambda }{2{\kern 1pt} s} )=\frac{2}{\lambda } e^{-\, \lambda /2{\kern 1pt} s} \sinh (\frac{\lambda }{2{\kern 1pt} s} )=\frac{1}{\lambda } \, (1-e^{-\, \lambda /s} )} \\ {\int _{0}^{\infty }e^{-\, s{\kern 1pt} t}  W_{1,2} (-\lambda {\kern 1pt} t)\, dt=\frac{1}{s} \, E_{1,2} (-\frac{\lambda }{s} )} \\ {\tau =\frac{\lambda }{s} } \\ {E_{1,2} (-\, \tau )=\, \frac{1}{\tau } (1-e^{-\, \tau } )} \end{array}
\end{equation}
 In general case, the Laplace transform can be expressed in terms of the incomplete gamma function $\gamma(a,z)$ [12]
\begin{equation} \label{GrindEQ__3_13_}
\begin{array}{l} {\int _{0}^{\infty }e^{-\, s{\kern 1pt} t}  W_{1,\beta \, +\, 1} (-\lambda {\kern 1pt} t)\, dt=\frac{1}{\lambda ^{\beta /2} } \, \int _{0}^{\infty }e^{-\, s{\kern 1pt} t}  t^{-\, \beta /2} J_{\beta } (2\sqrt{\lambda {\kern 1pt} t} )\, dt=} \\ {\frac{e^{i{\kern 1pt} \pi {\kern 1pt} \beta } s^{\beta \, -\, 1} }{\lambda ^{\beta } \Gamma (\beta )} \, e^{-\, \lambda /s} \gamma (\beta ,\frac{\lambda }{s} \, e^{-\, i{\kern 1pt} \pi {\kern 1pt} \beta } )\quad ;\quad Res>0} \end{array}
\end{equation}
and therefore
\begin{equation} \label{GrindEQ__3_14_}
\begin{array}{l} {L\left\{W_{1,\beta \, +\, 1} (-\lambda {\kern 1pt} t)\right\}=\frac{1}{s} \, E_{1,\beta \, +\, 1} (-\frac{\lambda }{s} )=\frac{e^{i{\kern 1pt} \pi {\kern 1pt} \beta } s^{\beta \, -\, 1} }{\lambda ^{\beta } \Gamma (\beta )} \, e^{-\, \lambda /s} \gamma (\beta ,\frac{\lambda }{s} \, e^{-\, i{\kern 1pt} \pi {\kern 1pt} \beta } )} \\ {z=\frac{\lambda }{s} } \\ {E_{1,\beta \, +\, 1} (-\, z)=\frac{e^{i{\kern 1pt} \pi {\kern 1pt} \beta } }{\Gamma (\beta )\, z^{\beta } } \, e^{-\, z} \gamma (\beta ,z\, e^{-\, i{\kern 1pt} \pi {\kern 1pt} \beta } )} \\ {} \end{array}
\end{equation}
For $\beta = 2$, we have $\mathrm{exp}(\pm 2i\pi)$ and
\begin{equation} \label{GrindEQ__3_15_}
E_{1,3} (-\, z)=\frac{1}{z^{2} } \, e^{-\, z} \gamma (2,z)=\frac{1}{z^{2} } \, e^{-\, z} \int _{0}^{z}e^{-{\kern 1pt} t} t\, dt
\end{equation}   

If $n$ is positive integer, then
\begin{equation} \label{GrindEQ__3_16_}
\begin{array}{l} {\gamma (n,z)=\Gamma (n){\kern 1pt} P(n,z)} \\ {P(n,z)=1-\left(1+z+\frac{z^{2} }{2!} +...\frac{z^{n\, -\, 1} }{(n-1)!} \right)\, e^{-\_ z} } \\ {\gamma (2,z)=1-(1+z){\kern 1pt} \, e^{-\, z} } \end{array}
\end{equation}
There are some equivalent expressions in the form given in [12]
\begin{equation} \label{GrindEQ__3_17_}
\begin{array}{l} {E_{1,\beta \, +\, 1} (-\, z)=\frac{e^{i{\kern 1pt} \pi {\kern 1pt} \beta } }{\Gamma (\beta )\, z^{\beta } } \, e^{-\, z} \gamma (\beta ,z\, e^{-\, i{\kern 1pt} \pi {\kern 1pt} \beta } )\, } \\ {E_{1,\beta \, +\, 1} (-\, z)=\frac{\sqrt{\pi } \, e^{-\, z/2} }{\Gamma (\beta +1)\, z^{(\beta \, +\, 1)/2} } \, M_{(1\, -\, \beta )/2,\beta /2} (z)} \\ {E_{1,\beta \, +\, 1} (-\, z)=\frac{1}{\Gamma (\beta +1)\, } \, _{1} F_{1} (1;\beta +1;-\, z)=\frac{1}{\Gamma (\beta )\, } \, \int _{0}^{1}e^{z{\kern 1pt} t}  (1-t)^{\beta \, -\, 1} dt} \end{array}
\end{equation}
For $\beta$ being positive integer $n$, the last equation links the Mittag-Leffler functions with the Kummer functions (see also Appendix A in [14] for other results).

For positive values of argument $t$ we have
\begin{equation} \label{GrindEQ__3_18_}
W_{1,\beta \, +\, 1} (t)=t^{\beta \, /2} I_{\beta \, } (2\sqrt{t} )
\end{equation}
and therefore
\begin{equation} \label{GrindEQ__3_19_}
\int _{0}^{\infty }e^{-\, s{\kern 1pt} t}  W_{1,\beta \, +\, 1} (\lambda {\kern 1pt} t)\, dt=\frac{1}{\lambda ^{\beta /2} } \, \int _{0}^{\infty }e^{-\, s{\kern 1pt} t}  t^{-\, \beta /2} I_{\beta } (2\sqrt{\lambda {\kern 1pt} t} )\, dt
\end{equation}
For $\beta = 0$, this gives [14]
\begin{equation} \label{GrindEQ__3_20_}
\begin{array}{l} {\int _{0}^{\infty }e^{-\, s{\kern 1pt} t}  W_{1,1} (\lambda {\kern 1pt} t)\, dt=\, \int _{0}^{\infty }e^{-\, s{\kern 1pt} t}  I_{0} (2\sqrt{\lambda {\kern 1pt} t} )\, dt=\frac{1}{s} \, e^{\lambda /s} \quad } \\ {Res>0} \end{array}
\end{equation}
but
\begin{equation} \label{GrindEQ__3_21_}
L\left\{W_{1,1} (\lambda {\kern 1pt} t)\right\}=\frac{1}{s} \, E_{1,1} (\frac{\lambda }{s} )
\end{equation}
and by comparing the expected result is reached
\begin{equation} \label{GrindEQ__3_22_}
\begin{array}{l} {E_{1,1} (\frac{\lambda }{s} )=e^{\lambda /s} } \\ {z=\frac{\lambda }{s} } \\ {E_{1,1} (z)=e^{\, z} } \end{array}
\end{equation}
If $\beta = 1$, then [12]
\begin{equation} \label{GrindEQ__3_23_}
\begin{array}{l} {\int _{0}^{\infty }e^{-\, s{\kern 1pt} t}  W_{1,2} (\lambda {\kern 1pt} t)\, dt=\, \frac{1}{\sqrt{\lambda } } \int _{0}^{\infty }e^{-\, s{\kern 1pt} t}  t^{-\, 1/2} I_{1} (2\sqrt{\lambda {\kern 1pt} t} )\, dt=} \\ {\frac{1}{\lambda } \, (e^{\lambda /s} -1)} \\ {\int _{0}^{\infty }e^{-\, s{\kern 1pt} t}  W_{1,2} (-\lambda {\kern 1pt} t)\, dt=\frac{1}{s} \, E_{1,2} (\frac{\lambda }{s} )} \\ {z=\frac{\lambda }{s} } \\ {E_{1,2} (\, z)=\, \frac{1}{z} (e^{\, z} -1)} \end{array}
\end{equation}
In general case [12]
\begin{equation} \label{GrindEQ__3_24_}
\begin{array}{l} {\int _{0}^{\infty }e^{-\, s{\kern 1pt} t}  W_{1,\beta \, +\, 1} (\lambda {\kern 1pt} t)\, dt=\frac{1}{\lambda ^{\beta /2} } \, \int _{0}^{\infty }e^{-\, s{\kern 1pt} t}  t^{-\, \beta /2} \, I_{\beta } (2\sqrt{\lambda {\kern 1pt} t} )\, dt=} \\ {\frac{e^{\lambda /s} s^{\beta \, -\, 1} }{\Gamma (\beta )\, \lambda ^{\beta } } \, \gamma (\beta ,\frac{\lambda }{s} )} \\ {L\left\{W_{1,\beta \, +\, 1} (\lambda {\kern 1pt} t)\, \right\}=\frac{1}{s} \, E_{1,\beta \, +\, 1} (\frac{\lambda }{s} )} \\ {z=\frac{\lambda }{s} } \\ {E_{1,\beta \, +\, 1} (z)=\frac{e^{z} }{z^{\beta } } \, \gamma (\beta ,z)} \end{array}
\end{equation}
where the incomplete gamma function can be expressed in terms of the Kummer function
\begin{equation} \label{GrindEQ__3_25_}
\gamma (\beta ,z)=\frac{z^{\beta } }{\beta } \, e^{-\, z} \, _{1} F_{1} (1;\beta +1;z)=\frac{z^{\beta } }{\beta } \, \, _{1} F_{1} (1;\beta +1;-\, z)
\end{equation}
From \eqref{GrindEQ__2_24_} and \eqref{GrindEQ__2_25_} it follows that
\begin{equation} \label{GrindEQ__3_26_}
E_{1,\beta \, +\, 1} (z)=\frac{e^{z} }{z^{\beta } } \, \gamma (\beta ,z)=\frac{1}{\beta } \, _{1} F_{1} (1;\beta +1;z)=\frac{e^{z} }{\beta } \, _{1} F_{1} (1;\beta +1;-\, z)
\end{equation}
Particular values of the incomplete gamma function of interest are
\begin{equation} \label{GrindEQ__3_27_}
\begin{array}{l} {\gamma (1,z)=(1-e^{-\, z} )} \\ {E_{1,2} (z)=\frac{1}{z} (e^{\, z} -1)} \end{array}
\end{equation}
and
\begin{equation} \label{GrindEQ__3_28_}
\begin{array}{l} {\gamma (1/2,z)=\sqrt{\pi } \, erf(\sqrt{z} )\, } \\ {E_{1,3/2} (z)=\sqrt{\frac{\pi }{z} } \, e^{z} erf(\sqrt{z} )\, } \end{array}
\end{equation}
In derivation explicit expressions for the Mittag-Leffler functions the recurrence relation of the incomplete gamma function
\begin{equation} \label{GrindEQ__3_29_}
\gamma (a+1,z)=a{\kern 1pt} \gamma (a,z)-z^{a} e^{-\, z}
\end{equation}
is very useful. For example, for $n = 1,2,3$ we have
\begin{equation} \label{GrindEQ__3_30_}
\begin{array}{l} {\gamma (1,z)=\frac{1}{z} (1-e^{-\, z} )} \\ {\gamma (2,z)={\kern 1pt} \gamma (1,z)-z{\kern 1pt} e^{-\, z} =\frac{1}{z} (1-e^{-\, z} )-z{\kern 1pt} e^{-\, z} } \\ {\gamma (3,z)={\kern 1pt} \gamma (2,z)-z{\kern 1pt} e^{-\, z} =\frac{1}{z} (1-e^{-\, z} )-2z{\kern 1pt} e^{-\, z} } \\ {\gamma (n+1,z)=n{\kern 1pt} {\kern 1pt} \gamma (n,z)-z{\kern 1pt} e^{-\, z} \quad ;\quad n=1,2,3,...} \end{array}
\end{equation}
and previously derived formulas in \eqref{GrindEQ__3_22_} and in \eqref{GrindEQ__3_23_} are reached.

\noindent For $n + 1/2$, from \eqref{GrindEQ__3_29_} it follows that
\begin{equation} \label{GrindEQ__3_31_}
\begin{array}{l} {\gamma (1/2,z)=\sqrt{\pi } \, erf(\sqrt{z} )} \\ {\gamma (3/2,z)=\sqrt{\pi } \, erf(\sqrt{z} )-z{\kern 1pt} e^{-\, z} \, } \\ {\gamma (5/2,z)=2{\kern 1pt} {\kern 1pt} [\sqrt{\pi } \, erf(\sqrt{z} )-{\kern 1pt} z{\kern 1pt} e^{-\, z} ]-z^{2} e^{-\, z} \, } \end{array}
\end{equation}
and immediately this gives by using \eqref{GrindEQ__3_24_}
\begin{equation} \label{GrindEQ__3_32_}
\begin{array}{l} {E_{1,3/2} (z)=\sqrt{\frac{\pi }{z} } \, e^{z} {\kern 1pt} erf(\sqrt{z} )} \\ {E_{1,5/2} (z)=\frac{e^{z} {\kern 1pt} }{z^{3/2} } \, [\sqrt{\pi } \, e^{z} erf(\sqrt{z} )-z\, e^{-\, z} ]\, } \\ {E_{1,7/2} (z)=\frac{e^{z} {\kern 1pt} }{z^{5/2} } \left\{\, 2\, [\sqrt{\pi } \, e^{z} erf(\sqrt{z} )-ze^{-\, z} {\kern 1pt} ]\, -z^{2} \, e^{-\, z} \right\}} \end{array}
\end{equation}

The Laplace transforms of the Mainardi functions are
\begin{equation} \label{GrindEQ__3_33_}
\begin{array}{l} {L\left\{\frac{1}{t} F_{\alpha } (\frac{\lambda }{t^{\alpha } } )\right\}=L\left\{\frac{\alpha {\kern 1pt} \lambda }{t^{\alpha \, +\, 1} } M_{\alpha } (\frac{\lambda }{t^{\alpha } } )\right\}=e^{-\, \lambda {\kern 1pt} s^{\alpha } } \quad } \\ {0<\alpha <1\quad ;\quad \lambda >0} \end{array}
\end{equation}
and
\begin{equation} \label{GrindEQ__3_34_}
\begin{array}{l} {L\left\{\frac{1}{\alpha } F_{\alpha } (\frac{\lambda }{t^{\alpha } } )\right\}=L\left\{\frac{\lambda }{t^{\alpha } } M_{\alpha } (\frac{\lambda }{t^{\alpha } } \right\}=\frac{\lambda }{s^{1\, -\, \alpha } } e^{-\, \lambda \, s^{\alpha } } \quad } \\ {0<\alpha <1\quad ;\quad \lambda >0\quad } \end{array}
\end{equation}
or the tern of the Wright function
\begin{equation} \label{GrindEQ__3_35_}
\begin{array}{l} {L\left\{\frac{1}{\alpha } W_{-\, \alpha ,0} (-\frac{\lambda }{t^{\alpha } } )\right\}=L\left\{\frac{\lambda }{t^{\alpha } } W_{-\, \alpha ,1\, -\, \alpha } (-\frac{\lambda }{t^{\alpha } } )\right\}=\frac{\lambda }{s^{1\, -\, \alpha } } e^{-\, \lambda \, s^{\alpha } } \quad } \\ {0<\alpha <1\quad ;\quad \lambda >0} \end{array}
\end{equation}
The inverse Laplace transforms are known only for $\alpha = 1/2$ and $\alpha = 1/3$.
\begin{equation} \label{GrindEQ__3_36_}
\begin{array}{l} {L\left\{\frac{1}{t} F_{1/2} (\frac{\lambda }{t^{1/2} } )\right\}=L\left\{\frac{{\kern 1pt} \lambda }{2{\kern 1pt} t^{3/2\, } } M_{1/2} (\frac{\lambda }{t^{1/2} } )\right\}=e^{-\, \lambda {\kern 1pt} s^{1/2} } } \\ {\lambda >0\quad } \\ {L^{-\, 1} \left\{e^{-\, \lambda {\kern 1pt} s^{1/2} } \right\}=\frac{\lambda {\kern 1pt} e^{-\, \lambda ^{2} /4{\kern 1pt} t} }{2\sqrt{\pi } \, t^{3/2} } } \\ {\frac{1}{t} F_{1/2} (\frac{\lambda }{t^{1/2} } )=\frac{{\kern 1pt} \lambda }{2{\kern 1pt} t^{3/2\, } } M_{1/2} (\frac{\lambda }{t^{1/2} } )=\frac{\lambda {\kern 1pt} e^{-\, \lambda ^{2} /4{\kern 1pt} t} }{2\sqrt{\pi } \, t^{3/2} } } \\ {F_{1/2} (\lambda {\kern 1pt} \tau )=\frac{{\kern 1pt} \lambda {\kern 1pt} \tau ^{2} }{2{\kern 1pt} } M_{1/2} (\lambda {\kern 1pt} \tau )=\frac{\lambda {\kern 1pt} \tau ^{2} e^{-\, \lambda ^{2} \tau ^{2} /4{\kern 1pt} } }{2\sqrt{\pi } } } \end{array}
\end{equation}
Multiplication by $t$ of the Mainardi function in \eqref{GrindEQ__3_36_} is equivalent to
\begin{equation} \label{GrindEQ__3_37_}
\begin{array}{l} {L\left\{F_{1/2} (\frac{\lambda }{t^{1/2} } )\right\}=L\left\{\frac{{\kern 1pt} \lambda }{2{\kern 1pt} t^{1/2} } M_{1/2} (\frac{\lambda }{t^{1/2} } )\right\}=-\frac{d}{ds} \left\{e^{-\, \lambda {\kern 1pt} s^{1/2} } \right\}=\frac{\lambda }{2\, s^{1/2} } e^{-\, \lambda {\kern 1pt} s^{1/2} } } \\ {L^{-\, 1} \left\{\frac{\lambda }{2\, s^{1/2} } e^{-\, \lambda {\kern 1pt} s^{1/2} } \right\}=\frac{\lambda \, e^{-\, \lambda ^{2} /4{\kern 1pt} t} }{2\sqrt{\pi {\kern 1pt} t} \, } } \\ {F_{1/2} (\frac{\lambda }{t^{1/2} } )=\frac{{\kern 1pt} \lambda }{2{\kern 1pt} t^{1/2\, } } M_{1/2} (\frac{\lambda }{t^{1/2} } )=\frac{\lambda e^{-\, \lambda ^{2} /4{\kern 1pt} t} }{2\sqrt{\pi {\kern 1pt} t} \, } } \\ {F_{1/2} (\lambda {\kern 1pt} \tau )=\frac{{\kern 1pt} \lambda {\kern 1pt} \tau }{2{\kern 1pt} } M_{1/2} (\lambda {\kern 1pt} \tau )=\frac{\lambda \tau \, e^{-\, \lambda ^{2} \tau /4{\kern 1pt} } }{2\sqrt{\pi } } } \end{array}
\end{equation}
The same results, but in terms of the Wright functions can be written as
\begin{equation} \label{GrindEQ__3_38_}
\begin{array}{l} {L\left\{2\; W_{-\, 1/2,0} (-\frac{\lambda }{t^{1/2} } )\right\}=L\left\{\frac{\lambda }{t^{1/2} } W_{-\, 1/2,1\, /2} (-\frac{\lambda }{t^{1/2} } )\right\}=\frac{\lambda }{s^{1\, /2} } e^{-\, \lambda \, s^{1/2} } } \\ {L^{-\, 1} \left\{\frac{\lambda }{s^{1\, /2} } e^{-\, \lambda \, s^{1/2} } \right\}=\frac{1}{\sqrt{\pi {\kern 1pt} t} } \, e^{-\, \lambda ^{2} /4{\kern 1pt} t} } \\ {2\; W_{-\, 1/2,0} (-\frac{\lambda }{t^{1/2} } )=\frac{\lambda }{t^{1/2} } W_{-\, 1/2,1\, /2} (-\frac{\lambda }{t^{1/2} } )=\frac{\lambda }{\sqrt{\pi {\kern 1pt} t} } \, e^{-\, \lambda ^{2} /4{\kern 1pt} t} } \\ {2\; W_{-\, 1/2,0} (-\lambda {\kern 1pt} \tau )=\lambda {\kern 1pt} \tau \, W_{-\, 1/2,1\, /2} (-\lambda {\kern 1pt} \tau )=\frac{\lambda {\kern 1pt} \tau }{\sqrt{\pi {\kern 1pt} } } \, e^{-\, \lambda ^{2} {\kern 1pt} \tau ^{2} /4{\kern 1pt} } } \end{array}
\end{equation}

In general case of the multiplication by $t^{n}$,  the differentiation of exponential functions can be expressed in terms of the Bessel functions [12]        
\begin{equation} \label{GrindEQ__3_39_}
\begin{array}{l} {L\left\{t^{n} F_{1/2} (\frac{\lambda }{t^{1/2} } )\right\}=L\left\{\frac{{\kern 1pt} \lambda \, {\kern 1pt} t^{n\, -\, 1/2} }{2{\kern 1pt} } M_{1/2} (\frac{\lambda }{t^{1/2} } )\right\}=(-1)^{n} \frac{d^{n} }{ds^{n} } \left\{e^{-\, \lambda {\kern 1pt} s^{1/2} } \right\}=} \\ {\frac{\lambda ^{n\, \, +\, \, 1/2} \, s^{(1\, -\, 2{\kern 1pt} n)/4} }{2^{n\, \, -\, \, 1/2} \, \sqrt{\pi } } \, K_{n\, \, -\, \, 1/2} \, (\lambda \, s^{1/2} )} \end{array}       \eqref{GrindEQ__4_39_}
\end{equation}

and therefore from \eqref{GrindEQ__3_39_} we have
\begin{equation} \label{GrindEQ__3_40_}
\begin{array}{l} {L\left\{2{\kern 1pt} {\kern 1pt} t^{n} \; W_{-\, 1/2,0} (-\frac{\lambda }{t^{1/2} } )\right\}=L\left\{\lambda \, t^{n\, -\, 1/2} {\kern 1pt} W_{-\, 1/2,1\, /2} (-\frac{\lambda }{t^{1/2} } )\right\}=} \\ {(-1)^{n} \lambda \frac{d^{n} }{d{\kern 1pt} s^{n} } \left\{\frac{e^{-\, \lambda \, s^{1/2} } }{s^{1\, /2} } \right\}=\frac{\lambda ^{n\, +\, 3/2} \, s^{-(2{\kern 1pt} n\, +1\, )/4} }{2^{n\, -\, 1/2} \sqrt{\pi } } \, K_{n\, +\, 1/2} (\lambda {\kern 1pt} s^{1/2} )} \end{array}
\end{equation}

If $\alpha = 1/3$, then [3,4]
\begin{equation} \label{GrindEQ__3_41_}
L\left\{\frac{1}{t} F_{1/3} (\frac{\lambda }{t^{1/3} } )\right\}=L\left\{\frac{{\kern 1pt} \lambda }{3{\kern 1pt} t^{4/3\, } } M_{1/3} (\frac{\lambda }{t^{1/3} } )\right\}=e^{-\, \lambda {\kern 1pt} s^{1/3} } \quad
\end{equation}
but using [15]
\begin{equation} \label{GrindEQ__3_42_}
L\left\{\frac{\lambda ^{3/2} }{3{\kern 1pt} \pi {\kern 1pt} t^{3/2} } \, K_{1/3} \left(\frac{2{\kern 1pt} \lambda ^{3/2} }{\sqrt{27{\kern 1pt} t} } \right)\right\}=e^{-\, \lambda {\kern 1pt} s^{1/3} }
\end{equation}
we have
\begin{equation} \label{GrindEQ__3_43_}
3{\kern 1pt} F_{1/3} (\frac{\lambda }{t^{1/3} } )=\frac{{\kern 1pt} \lambda }{t^{1/3\, } } M_{1/3} (\frac{\lambda }{t^{1/3} } )=\frac{\lambda ^{3/2} }{\pi {\kern 1pt} t^{1/2} } \, K_{1/3} (\frac{2{\kern 1pt} \lambda ^{3/2} }{\sqrt{27t} } )\quad
\end{equation}
The same result is available from [3,4]
\begin{equation} \label{GrindEQ__3_44_}
L\left\{3{\kern 1pt} F_{1/3} (\frac{\lambda }{t^{1/3} } )\right\}=L\left\{\frac{\lambda }{t^{1/3} } M_{1/3} (\frac{\lambda }{t^{1/3} } \right\}=\frac{\lambda }{s^{2/3\, } } e^{-\, \lambda \, s^{\alpha } } \quad
\end{equation}
and [15]
\begin{equation} \label{GrindEQ__3_45_}
L\left\{\frac{\lambda ^{3/2} }{{\kern 1pt} \pi {\kern 1pt} t^{1/2} } \, K_{1/3} \left(\frac{2{\kern 1pt} \lambda ^{3/2} }{\sqrt{27{\kern 1pt} t} } \right)\right\}=\frac{\lambda {\kern 1pt} e^{-\, \lambda {\kern 1pt} s^{1/3} } }{s^{2/3} }
\end{equation}
In terms of the Wright functions it can be expressed by
\begin{equation} \label{GrindEQ__3_46_}
3{\kern 1pt} {\kern 1pt} W_{-\, 1/3,0} (-\frac{\lambda }{t^{1/3} } )=\frac{\lambda }{t^{1/3} } W_{-\, 1/3,2/3} (-\frac{\lambda }{t^{1/3} } )=\frac{\lambda ^{3/2} }{\pi {\kern 1pt} t^{1/2} } K_{1/3} (\frac{2{\kern 1pt} \lambda ^{3/2} }{\sqrt{27{\kern 1pt} t} } )
\end{equation}

\section{Functional limits associated with the Wright \\ functions.}

In 1969 Lamborn, see see [8-11] and [18-20]
 proposed the following delta sequence for representation of the shifted Dirac delta function
\begin{equation} \label{GrindEQ__4_1_}
\delta (x-1)={\mathop{\lim }\limits_{\nu \, \to \, \infty }} \left[\nu {\kern 1pt} J_{\nu } (\nu x{\kern 1pt} )\right] 
\end{equation}
As it was demonstrated over the 2000-2008 period by Apelblat [16-18], this delta sequence is useful for evaluation of the asymptotic relations, limits of series, integrals and integral representations of elementary and special functions.

If the Lamborn expression is multiplied by a function \textit{f}(\textit{tx}) and integrated from zero to infinity with respect to variable \textit{x} we have
\begin{equation} \label{GrindEQ__4_2_}
f(t)=\int _{0}^{\infty }f(t{\kern 1pt} x )\delta (x-1){\kern 1pt} {\kern 1pt} dx={\mathop{\lim }\limits_{\nu \, \to \, \infty }} \left[\nu \, \int _{0}^{\infty }f(t{\kern 1pt} x)J_{\nu } (\nu {\kern 1pt} x){\kern 1pt} {\kern 1pt} dx \right]
\end{equation}
In such way, the function $f(t)$ is represented by the asymptotic limit of the infinite integral of product of $f(tx)$ and the Bessel function $J_{\nu}(\nu x)$. If the right-hand integral in \eqref{GrindEQ__4_2_} can be evaluated in the closed form, then the limit can be regarded as the generalization of the l'Hospital's rule.
\begin{equation} \label{GrindEQ__4_3_}
\begin{array}{l} {f(t)={\mathop{\lim }\limits_{\nu \, \to \, \infty }} [\nu {\kern 1pt} {\kern 1pt} \Phi (t,\nu )]} \\ {\Phi (t,\nu )=\int _{0}^{\infty }f(t{\kern 1pt} x)J_{\nu }  (\nu {\kern 1pt} x){\kern 1pt} {\kern 1pt} dx} \end{array}
\end{equation}
For the Wright function treated as the generalized the Bessel function
\begin{equation} \label{GrindEQ__4_4_}
f(t)=W_{1,\beta \, +\, 1} (-\frac{t^{2} }{4} )=\left(\frac{2}{t} \right)^{\beta } J_{\beta } (t)
\end{equation}
it follows from ((4.1) and \eqref{GrindEQ__4_4_} that
\begin{equation} \label{GrindEQ__4_5_}
\begin{array}{l} {f(t)={\mathop{\lim }\limits_{\nu \, \to \, \infty }} \left\{\nu \, \int _{0}^{\infty }J_{\nu }  (\nu {\kern 1pt} x)\left(\frac{2}{t{\kern 1pt} x} \right)^{\chi } J_{\beta } (t{\kern 1pt} x)\, dx\, \right\}=W_{1,\beta \, +\, 1} \left(-\frac{t^{2} }{4} \right)} \\ {\Phi (t,\nu ,\beta )=\int _{0}^{\infty }x^{-\, \beta } J_{\nu }  (\nu {\kern 1pt} x)\, J_{\beta } (t{\kern 1pt} x)\, dx} \end{array}
\end{equation}
However, the infinite integral in \eqref{GrindEQ__4_5_} is known [19]
\begin{equation} \label{GrindEQ__4_6_}
\Phi (t,\nu ,\mu )=\int _{0}^{\infty }x^{-\, \beta } J_{\nu }  (\nu {\kern 1pt} x)\, J_{\beta } (t{\kern 1pt} x)\, dx=\left(\frac{t}{2} \right)^{\beta } \frac{1}{\nu } \, _{2} F_{1} \left(\frac{\nu +1}{2} ,\frac{1-\nu }{2} ;\beta +1;\frac{t^{2} }{\nu ^{2} } \right)
\end{equation}
and therefore the Wright function is represented by the following limit
\begin{equation} \label{GrindEQ__4_7_}
\begin{array}{l} {W_{1,\beta \, +\, 1} \left(-\frac{t^{2} }{4} \right)={\mathop{\lim }\limits_{\nu \, \to \, \infty }} \, \left\{{}_{2} F_{1} \left(\frac{\nu +1}{2} ,\frac{1-\nu }{2} ;\beta +1\, ;\frac{t^{2} }{\nu ^{2} } \right)\right\}\, } \\ {Re(\nu +1)>0\quad ;\quad Re\beta >-1\quad ;\quad 0<t<\nu } \end{array}
\end{equation}
or in the equivalent form
\begin{equation} \label{GrindEQ__4_8_}
W_{1,\beta \, +\, 1} \left(-\, x\right)={\mathop{\lim }\limits_{\nu \, \to \, \infty }} \, \left\{{}_{2} F_{1} \left(\frac{\nu +1}{2} ,\frac{1-\nu }{2} ;\beta +1\, ;\frac{4{\kern 1pt} x}{\nu ^{2} } \right)\right\}\,
\end{equation}
For $\beta=0$, we have
\begin{equation} \label{GrindEQ__4_9_}
W_{1,\, 1} \left(-\frac{t^{2} }{4} \right)=J_{0} (t)={\mathop{\lim }\limits_{\nu \, \to \, \infty }} \, \left\{{}_{2} F_{1} \left(\frac{\nu +1}{2} ,\frac{1-\nu }{2} ;1\, ;\frac{t^{2} }{\nu ^{2} } \right)\right\}\,
\end{equation}
and for $\beta = \pm 1/2$.
\begin{equation} \label{GrindEQ__4_10_}
\begin{array}{l} {W_{1,\, 3/2} \left(-\frac{t^{2} }{4} \right)=\sqrt{\frac{2}{t} } J_{1/2} (t)={\mathop{\lim }\limits_{\nu \, \to \, \infty }} \, \left\{{}_{2} F_{1} \left(\frac{\nu +1}{2} ,\frac{1-\nu }{2} ;\frac{3}{2} ;\frac{t^{2} }{\nu ^{2} } \right)\right\}=\frac{2{\kern 1pt} {\kern 1pt} \sin t}{\sqrt{\pi } \_ \, t} \, } \\ {W_{1,\, 1/2} \left(-\frac{t^{2} }{4} \right)=\sqrt{\frac{t}{2} } J_{-\, 1/2} (t)={\mathop{\lim }\limits_{\nu \, \to \, \infty }} \, \left\{{}_{2} F_{1} \left(\frac{\nu +1}{2} ,\frac{1-\nu }{2} ;\frac{1}{2} \, ;\frac{t^{2} }{\nu ^{2} } \right)\right\}=-\frac{\cos t}{\sqrt{\pi } \_ \, } \, } \end{array}
\end{equation}
Hypergeometric functions \eqref{GrindEQ__4_10_} are known in different form [20]
\begin{equation} \label{GrindEQ__4_11_}
\begin{array}{l} {{}_{2} F_{1} \left(a,1-a;\frac{3}{2} ;(\sin z)^{2} \right)=\frac{\sin [(2{\kern 1pt} a-1)z]}{(2{\kern 1pt} a-1)\sin z} } \\ {{}_{2} F_{1} \left(a,1-a;\frac{1}{2} ;(\sin z)^{2} \right)=\frac{\cos [(2{\kern 1pt} a-1)z]}{\cos z} } \\ {a=\frac{\nu +1}{2} \quad ;\quad \sin z=\frac{t}{\nu } } \end{array}
\end{equation}

If the delta sequence in \eqref{GrindEQ__4_1_} is used together with integral transforms having different kernels $T$, the we have [17,18]
\begin{equation} \label{GrindEQ__4_12_}
\begin{array}{l} {{\mathop{\lim }\limits_{\nu \, \to \, \infty }} [\nu \, \int _{0}^{\infty }f(\xi ,\lambda ){\kern 1pt} T\{ J_{\nu } (\nu {\kern 1pt} x),\xi \} \, d\xi ] =} \\ {{\mathop{\lim }\limits_{\nu \, \to \, \infty }} [\nu \, \int _{0}^{\infty }J_{\nu } (\nu {\kern 1pt} x){\kern 1pt} T\{ f(\xi ,\lambda ),x\} \, dx] =T(1,\lambda )} \end{array}
\end{equation}
In the case of the Laplace transformation, \eqref{GrindEQ__4_12_} can be written in the following way
\begin{equation} \label{GrindEQ__4_13_}
\begin{array}{l} {\int _{0}^{\infty }e^{-\, \xi {\kern 1pt} x} J_{\nu } (\nu {\kern 1pt} \xi )\, d\xi =\frac{\nu ^{\nu } }{\sqrt{\nu ^{2} +\xi ^{2} } \, [\xi +\sqrt{\nu ^{2} +\xi ^{2} } ]^{\nu } }  } \\ {{\mathop{\lim }\limits_{\nu \, \to \infty \, }} \left\{\nu ^{\nu \, +\, 1} \, \int _{0}^{\infty }\frac{f(\xi ,\lambda )}{\sqrt{\nu ^{2} +\xi ^{2} } \, [\xi +\sqrt{\nu ^{2} +\xi ^{2} } ]^{\nu } } \, d\xi  \right\}=L(1,\lambda )} \end{array}
\end{equation}
Introducing \eqref{GrindEQ__3_1_} into \eqref{GrindEQ__4_13_} we have
\begin{equation} \label{GrindEQ__4_14_}
{\mathop{\lim }\limits_{\nu \, \to \infty \, }} \left\{\nu ^{\nu \, +\, 1} \, \int _{0}^{\infty }\frac{W_{\alpha ,\beta } (\lambda {\kern 1pt} \xi )}{\sqrt{\nu ^{2} +\xi ^{2} } \, [\xi +\sqrt{\nu ^{2} +\xi ^{2} } ]^{\nu } } \, d\xi  \right\}=\left. \frac{1}{s} \, E_{\alpha ,\beta } (\frac{\lambda }{s} )\right|_{s\, =\, 1} \, =E_{\alpha ,\beta } (\lambda )
\end{equation}
The same operation performed with \eqref{GrindEQ__3_2_} gives
\begin{equation} \label{GrindEQ__4_15_}
\begin{array}{l} {{\mathop{\lim }\limits_{\nu \, \to \infty \, }} \left\{\nu ^{\nu \, +\, 1} \, \int _{0}^{\infty }\frac{e^{\rho {\kern 1pt} \xi } \, W_{\alpha ,\beta } (\lambda {\kern 1pt} \xi )}{\sqrt{\nu ^{2} +\xi ^{2} } \, [\xi +\sqrt{\nu ^{2} +\xi ^{2} } ]^{\nu } } \, d\xi  \right\}=\frac{1}{1-\rho } E_{\alpha ,} {}_{\beta } (\frac{\lambda }{1-\rho } )} \\ {{\mathop{\lim }\limits_{\nu \, \to \infty \, }} \left\{\nu ^{\nu \, +\, 1} \, \int _{0}^{\infty }\frac{e^{-\, \rho {\kern 1pt} \xi } \, W_{\alpha ,\beta } (\lambda {\kern 1pt} \xi )}{\sqrt{\nu ^{2} +\xi ^{2} } \, [\xi +\sqrt{\nu ^{2} +\xi ^{2} } ]^{\nu } } \, d\xi  \right\}=\frac{1}{1+\rho } E_{\alpha ,} {}_{\beta } (\frac{\lambda }{1+\rho } )} \end{array}
\end{equation}
and using \eqref{GrindEQ__3_6_}
\begin{equation} \label{GrindEQ__4_16_}
\begin{array}{l} {{\mathop{\lim }\limits_{\nu \, \to \infty \, }} \left\{\nu ^{\nu \, +\, 1} \, \int _{0}^{\infty }\frac{\xi \, W_{\alpha ,\beta } (\lambda {\kern 1pt} \xi )}{\sqrt{\nu ^{2} +\xi ^{2} } \, [\xi +\sqrt{\nu ^{2} +\xi ^{2} } ]^{\nu } } \, d\xi  \right\}=} \\ {\frac{1}{\alpha {\kern 1pt} \lambda } [(\alpha {\kern 1pt} \lambda -\beta +1)E_{\alpha ,\beta } (\lambda )+E_{\alpha ,\beta \, -\, 1} (\lambda )]} \end{array}
\end{equation}
The Laplace transform in \eqref{GrindEQ__3_24_} leads to
\begin{equation} \label{GrindEQ__4_17_}
{\mathop{\lim }\limits_{\nu \, \to \infty \, }} \left\{\nu ^{\nu \, +\, 1} \, \int _{0}^{\infty }\frac{W_{1,\beta \, +\, 1} (\lambda {\kern 1pt} \xi )}{\sqrt{\nu ^{2} +\xi ^{2} } \, [\xi +\sqrt{\nu ^{2} +\xi ^{2} } ]^{\nu } } \, d\xi  \right\}=\frac{e^{-\, \lambda } }{\lambda ^{\beta } \Gamma (\beta )} \gamma (\beta ,\lambda )\, =E_{1,\beta \, +\, 1} (\lambda )
\end{equation}
For integer values of parameters $\beta$ in \eqref{GrindEQ__4_17_}, the limits of the Wright functions can be represented by simple expressions. For $\beta = 1$, we have
\begin{equation} \label{GrindEQ__4_18_}
\begin{array}{l} {{\mathop{\lim }\limits_{\nu \, \to \infty \, }} \left\{\nu ^{\nu \, +\, 1} \, \int _{0}^{\infty }\frac{W_{1,1} (\lambda {\kern 1pt} \xi )}{\sqrt{\nu ^{2} +\xi ^{2} } \, [\xi +\sqrt{\nu ^{2} +\xi ^{2} } ]^{\nu } } \, d\xi  \right\}=e^{\lambda } } \\ {{\mathop{\lim }\limits_{\nu \, \to \infty \, }} \left\{\nu ^{\nu \, +\, 1} \, \int _{0}^{\infty }\frac{W_{1,1} (-\, \lambda {\kern 1pt} \xi )}{\sqrt{\nu ^{2} +\xi ^{2} } \, [\xi +\sqrt{\nu ^{2} +\xi ^{2} } ]^{\nu } } \, d\xi  \right\}=e^{-\, \lambda } } \end{array}
\end{equation}
and
\begin{equation} \label{GrindEQ__4_19_}
{\mathop{\lim }\limits_{\nu \, \to \infty \, }} \left\{\nu ^{\nu \, +\, 1} \, \int _{0}^{\infty }\frac{W_{1,1} (-\frac{\lambda ^{2} {\kern 1pt} \xi ^{2} }{4} )}{\sqrt{\nu ^{2} +\xi ^{2} } \, [\xi +\sqrt{\nu ^{2} +\xi ^{2} } ]^{\nu } } \, d\xi  \right\}=\frac{1}{\sqrt{1+\lambda ^{2} } }
\end{equation}
For $\beta = 1$, the corresponding limits are
\begin{equation} \label{GrindEQ__4_20_}
{\mathop{\lim }\limits_{\nu \, \to \infty \, }} \left\{\nu ^{\nu \, +\, 1} \, \int _{0}^{\infty }\frac{W_{1,2} (-\, \lambda {\kern 1pt} \xi )}{\sqrt{\nu ^{2} +\xi ^{2} } \, [\xi +\sqrt{\nu ^{2} +\xi ^{2} } ]^{\nu } } \, d\xi  \right\}=\frac{2}{\lambda } e^{-\, \lambda } \sinh (\frac{\lambda }{2} )
\end{equation}
\begin{equation} \label{GrindEQ__4_21_}
{\mathop{\lim }\limits_{\nu \, \to \infty \, }} \left\{\nu ^{\nu \, +\, 1} \, \int _{0}^{\infty }\frac{W_{1,2} (\lambda {\kern 1pt} \xi )}{\sqrt{\nu ^{2} +\xi ^{2} } \, [\xi +\sqrt{\nu ^{2} +\xi ^{2} } ]^{\nu } } \, d\xi  \right\}=\frac{1}{\lambda } (e^{\lambda } -1)
\end{equation}
\begin{equation} \label{GrindEQ__4_22_}
{\mathop{\lim }\limits_{\nu \, \to \infty \, }} \left\{\nu ^{\nu \, +\, 1} \, \int _{0}^{\infty }\frac{W_{1,2} (-\frac{\lambda ^{2} {\kern 1pt} \xi ^{2} }{4} )}{\sqrt{\nu ^{2} +\xi ^{2} } \, [\xi +\sqrt{\nu ^{2} +\xi ^{2} } ]^{\nu } } \, d\xi  \right\}=\frac{2}{1+\sqrt{1+\lambda ^{2} } }
\end{equation}
Similarly, for $\beta  = 3$,  the functional limit is
\begin{equation} \label{GrindEQ__4_23_}
\begin{array}{l} {{\mathop{\lim }\limits_{\nu \, \to \infty \, }} \left\{\nu ^{\nu \, +\, 1} \, \int _{0}^{\infty }\frac{W_{1,3} (-\frac{\lambda ^{2} {\kern 1pt} \xi ^{2} }{4} )}{\sqrt{\nu ^{2} +\xi ^{2} } \, [\xi +\sqrt{\nu ^{2} +\xi ^{2} } ]^{\nu } } \, d\xi  \right\}=} \\ {\frac{1}{\lambda } \left\{\frac{\lambda }{1+\sqrt{1+\lambda ^{2} } } +\frac{1}{3} \left[\frac{\lambda }{1+\sqrt{1+\lambda ^{2} } } \right]^{3} \right\}} \end{array}
\end{equation}
The Laplace transforms of the Mainardi functions from \eqref{GrindEQ__3_39_} and \eqref{GrindEQ__3_40_} are
\begin{equation} \label{GrindEQ__4_24_}
\begin{array}{l} {{\mathop{\lim }\limits_{\nu \, \to \infty \, }} \left\{\nu ^{\nu \, +\, 1} \, \int _{0}^{\infty }\frac{\xi ^{n} F_{1/2} (\frac{\lambda {\kern 1pt} }{\xi ^{1/2} } )}{\sqrt{\nu ^{2} +\xi ^{2} } \, [\xi +\sqrt{\nu ^{2} +\xi ^{2} } ]^{\nu } } \, d\xi  \right\}=\frac{\lambda ^{n\, +\, 1} }{2^{n\, -\, 1/2} \sqrt{\pi } } \, K_{n\, -\, 1/2} (\lambda )} \\ {{\mathop{\lim }\limits_{\nu \, \to \infty \, }} \left\{\nu ^{\nu \, +\, 1} \, \int _{0}^{\infty }\frac{\xi ^{n\, -\, 1/2} M_{1/2} (\frac{\lambda {\kern 1pt} }{\xi ^{1/2} } )}{\sqrt{\nu ^{2} +\xi ^{2} } \, [\xi +\sqrt{\nu ^{2} +\xi ^{2} } ]^{\nu } } \, d\xi  \right\}=\frac{\lambda ^{n\, -\, 1/2} }{2^{n\, +\, 1/2} \sqrt{\pi } } \, K_{n\, -\, 1/2} (\lambda )} \end{array}
\end{equation}
and from \eqref{GrindEQ__3_44_} we have
\begin{equation} \label{GrindEQ__4_25_}
{\mathop{\lim }\limits_{\nu \, \to \infty \, }} \left\{\nu ^{\nu \, +\, 1} \, \int _{0}^{\infty }\frac{F_{1/3} (\frac{\lambda {\kern 1pt} }{\xi ^{1/3} } )}{\sqrt{\nu ^{2} +\xi ^{2} } \, [\xi +\sqrt{\nu ^{2} +\xi ^{2} } ]^{\nu } } \, d\xi  \right\}=\frac{\lambda }{3} \, e^{-\, \lambda }
\end{equation}

\section{Conclusions}

 Parameters of the Wright functions were treated as variables and derivatives with respect to them were derived and discussed. These derivatives are expressible in terms of infinite power series with quotients of digamma and gamma functions in their coefficients. The functional form of these series resembles those which were derived for the Mittag-Leffler functions. Only in few cases, it was possible to obtain the sums of these series in a closed form. The differentiation operation when the Wright functions are treated as the generalized Bessel functions leads to the Bessel functions and their derivatives with respect to the order. Simple operations with the Laplace transforms of the Wright functions of the first kind give explicit forms of the Mittag-Leffler functions. Applying the shifted Dirac delta function, permits to derive functional limits by using the Laplace transforms of the Wright  functions
 
 Finally we would like to draw attention of the interested readers to the recent papers [23], [24], [25]  where some noteworthy  applications of the Wright functions of the first and of the second kind are discussed.
 {
 Relevant applications can be expected in the field of special functions in fractional calculus for which we address the interested readers to its extensive literature, see e.g.  the papers [26]-[32].}

\textbf{}

\noindent

\noindent

\noindent \textbf{Acknowledgements}

\noindent 
The
work of F. Mainardi has been carried out in the framework of the activities
of the National Group of Mathematical Physics (INdAM-GNFM) .
Both authors are grateful 
 to Associate Professor Juan Luis Gonzalez-Santander Martinez, Department of Mathematics, Universidad de Oviedo, Oviedo, Spain for his help with using MATHEMATICA program and with editing in LaTeX.

\section*{Appendix: Differentiation of the Wright functions \\
with respect to parameters  versus Bessel functions}

Initially, the Wright functions (of the first kind) were treated as the generalized Bessel functions because for parameters $\alpha = 1$ and $\beta + 1$ they become
\begin{equation} \label{GrindEQ__2_18_}
\begin{array}{l} {W_{1,\beta \, +\, 1} (-\frac{t^{2} }{4} )=\left(\frac{2}{t} \right)^{\beta } \, J_{\beta } (t)} \\ {W_{1,\beta \, +\, 1} (\frac{t^{2} }{4} )=\left(\frac{2}{t} \right)^{\beta } \, I_{\beta } (t)} \end{array}
\end{equation}
Differentiation of the Wright functions in \eqref{GrindEQ__2_18_} with respect to parameter $\beta$ gives
\begin{equation} \label{GrindEQ__2_19_}
\begin{array}{l} {\frac{\partial \left(W_{1,\beta \, +\, 1} (-\frac{t^{2} }{4} )\right)}{\partial \beta } =\left(\frac{2}{t} \right)^{\beta } \, \left[\ln \left(\frac{2}{t} \right)J_{\beta } (t)+\frac{\partial J_{\beta } (t)}{\partial \beta } \right]} \\ {\frac{\partial \left(W_{1,\beta \, +\, 1} (\frac{t^{2} }{4} )\right)}{\partial \beta } =\left(\frac{2}{t} \right)^{\beta } \, \left[\ln \left(\frac{2}{t} \right)I_{\beta } (t)+\, \frac{\partial I_{\beta } (t)}{\partial \beta } \right]} \end{array}
\end{equation}
However, differentiation of the Bessel functions with respect to the order can be expressed by [13]
\begin{equation} \label{GrindEQ__2_20_}
\begin{array}{l} {\frac{\partial J_{\beta } (t)}{\partial \beta } =\pi \, \beta \, \int _{0}^{\pi /2}\tan \theta \, \, Y_{0}  \left(t\, [\sin \theta ]^{2} \right)\, J_{\beta } \left(t\, [\cos \theta ]^{2} \right)\, d\theta } \\ {\frac{\partial I_{\beta } (t)}{\partial \beta } =-\, 2\, \beta \, \int _{0}^{\pi /2}\tan \theta \, \, K_{0}  \left(t\, [\sin \theta ]^{2} \right)\, I_{\beta } \left(t\, [\cos \theta ]^{2} \right)\, d\theta } \\ {Re\beta >0} \end{array}
\end{equation}
In particular cases, differentiation with respect to the parameter $\beta$ can be explicitly expressed [8]
\begin{equation} \label{GrindEQ__2_21_}
\begin{array}{l} {\left(\frac{\partial J_{\beta } (t)}{\partial \beta } \right)_{\beta \, =\, 0} =\frac{\pi }{2} \, Y_{0} (t)} \\ {\left(\frac{\partial I_{\beta } (t)}{\partial \beta } \right)_{\beta \, =\, 0} =-\, K_{0} (t)} \end{array}
\end{equation}
and therefore
\begin{equation} \label{GrindEQ__2_22_}
\begin{array}{l} {\left. \left(\frac{\partial {\kern 1pt} {\kern 1pt} W_{1,\beta \, +\, 1} (-\frac{t^{2} }{4} )}{\partial \beta } \right)\right|_{\beta \, =\, 0} =-\ln \left(\frac{t}{2} \right)J_{0} (t)+\frac{\pi }{2} Y_{0} (t)} \\ {\left. \left(\frac{\partial {\kern 1pt} {\kern 1pt} W_{1,\beta \, +\, 1} (\frac{t^{2} }{4} )}{\partial \beta } \right)\right|_{\beta \, =\, 0} =-\ln \left(\frac{t}{2} \right)I_{0} (t)-K_{0} (t)} \end{array}
\end{equation}
For $\beta = 1$ we have
\begin{equation} \label{GrindEQ__2_23_}
\begin{array}{l} {\left(\frac{\partial J_{\beta } (t)}{\partial \beta } \right)_{\beta \, =\, 1} =\frac{J_{0} (t)}{t} +\frac{\pi }{2} \, Y_{1} (t)} \\ {\left(\frac{\partial I_{\beta } (t)}{\partial \beta } \right)_{\beta \, =\, 1} =K_{1} (t)-\frac{I_{0} (t)}{t} } \end{array}
\end{equation}
which gives
\begin{equation} \label{GrindEQ__2_24_}
\begin{array}{l} {\left. \left(\frac{\partial {\kern 1pt} {\kern 1pt} W_{1,\beta \, +\, 1} (-\frac{t^{2} }{4} )}{\partial \beta } \right)\right|_{\beta \, =\, 1} =\left(\frac{2}{t} \right)\left[-\ln \left(\frac{t}{2} \right)J_{1} (t)-\frac{J_{0} (t)}{t} +\frac{\pi }{2} Y_{1} (t)\right]} \\ {\left. \left(\frac{\partial {\kern 1pt} {\kern 1pt} W_{1,\beta \, +\, 1} (\frac{t^{2} }{4} )}{\partial \beta } \right)\right|_{\beta \, =\, 1} =\left(\frac{2}{t} \right)\left[-\ln \left(\frac{t}{2} \right)I_{1} (t)+K_{1} (t)-\frac{I_{0} (t)}{t} \right]} \end{array}
\end{equation}
Derivatives for $beta = 1/2$ are
\begin{equation} \label{GrindEQ__2_25_}
\begin{array}{l} {\left(\frac{\partial J_{\beta } (t)}{\partial \beta } \right)_{\beta \, =\, 1/2} =\sqrt{\frac{2}{\pi \, t} } \, \left[\sin t\, Ci(2t)-\cos t\, Si(2t)\right]} \\ {\left(\frac{\partial I_{\beta } (t)}{\partial \beta } \right)_{\beta \, =\, 1/2} =\sqrt{\frac{1}{2\, \pi \, t} } \, \left[e^{t} Ei(-2t)-e^{-\, t} Ei(2t)\right]} \\ {J_{1/2} (t)=\sqrt{\frac{2}{\pi \, t} } \, \sin t\quad ;\quad I_{1/2} (t)=\sqrt{\frac{2}{\pi \, t} } \, \sinh t} \end{array}
\end{equation}
which leads to
\begin{equation} \label{GrindEQ__2_26_}
\begin{array}{l} {\left. \left(\frac{\partial {\kern 1pt} {\kern 1pt} W_{1,\beta \, +\, 1} (-\frac{t^{2} }{4} )}{\partial \beta } \right)\right|_{\beta =\, 1/2} =\frac{2}{\sqrt{\pi } \, t} \left[-\ln \left(\frac{t}{2} \right)\sin t+\sin t\, Ci(2t)-\cos t\, Si(2t)\right]} \\ {\left. \left(\frac{\partial {\kern 1pt} {\kern 1pt} W_{1,\beta \, +\, 1} (\frac{t^{2} }{4} )}{\partial \beta } \right)\right|_{\beta \, =\, 1/2} =\frac{2}{\sqrt{\pi } \, t} \left[-\ln \left(\frac{t}{2} \right)\sinh t+\frac{1}{2} \left(e^{t} Ei(-\, 2t)-e^{-\, t} Ei(2t)\right)\right]} \end{array}
\end{equation}
If variable is changed to $t = 2x^{1/2}$, these results can be equivalently written in different form, for example \eqref{GrindEQ__2_22_} is
\begin{equation} \label{GrindEQ__2_27_}
\begin{array}{l} {\left. \left(\frac{\partial {\kern 1pt} {\kern 1pt} W_{1,\beta \, +\, 1} (-x)}{\partial \beta } \right)\right|_{\beta \, =\, 0} =-\ln \sqrt{x} \, J_{0} (2\sqrt{x} )+\frac{\pi }{2} Y_{0} (2\sqrt{x} )} \\ {\left. \left(\frac{\partial {\kern 1pt} {\kern 1pt} W_{1,\beta \, +\, 1} (x)}{\partial \beta } \right)\right|_{\beta \, =\, 0} =-\ln \sqrt{x} \, I_{0} (2\sqrt{x} )-K_{0} (2\sqrt{x} )} \end{array}
\end{equation}

\noindent
\textbf{References}

\noindent [1] Wright, E.M., On the coefficients of power series having exponential  singularities. J. London Math. Soc., 1933, 8, 71-79.

\noindent [2] Wright, E.M., The generalized Bessel function of order greater than one. Quart. J.  Math. Oxford, 1940, 11, 36-48.

\noindent [3] Gorenflo, R., Luchko, Yu., Mainardi, F. Analytical properties and applications of  the Wright functions. 
Fract. Calc. Appl. Anal., 1999, 2, 383-414.

\noindent [4] Mainardi, F. Fractional Calculus and Waves in Linear Viscoelasticity, 2-nd edition, World Scientific, Singapore , 2022.
[First edition, Imperial   College Press, London, 2010].

\noindent [5] Kiryakova, V. A guide to special functions in fractional calculus. A survey,
Mathematics, 2021, 9, No 1, Art ID 106, 40 pp. DOI: 10.3390/math9010106.

\noindent [6] Srivastava, H.M. 
A survey of some recent developments on higher transcendental functions of analytic number theory and applied mathematics.
Symmetry, 2021, 13, 2294/1--22.   DOI: 10.3390/sym1312294
  
\noindent [7] Gorenflo, R., Kilbas, A,  Mainardi, F,, Rogosin, S.
Mittag-Leffler Functions, Related Topics and Applications,
2-nd edition, Springer, Heidelberg 2020. [First Edition 2014]
 
\noindent [8] Apelblat, A. Bessel and Related Functions. Mathematical Operations with  Respect to the Order. Vol 1.
 Theoretical Aspects. Walter de Gruyter GmbH,  Berlin, 2020.

\noindent [9] Apelblat, A. Bessel and Related Functions. Mathematical Operations with  Respect to the Order. Vol 2. Numerical Results. 
Walter de Gruyter GmbH,  Berlin, 2020.

\noindent [10] Apelblat, A., Differentiation of the Mittag-Leffler functions with respect to  parameters in the Laplace transform approach. 
Mathematics (MDPI), 2020, Vol 8,  Art.   657/1--22
DOI: 10.3390/math8050657


\noindent [11] Lamborn, B.N.  An expression for the Dirac delta function,
SIAM Rev. 1969,  11 (4), 603.
 
\noindent [12] Brychkov, Y.A. Handbook of Special Functions. Derivatives, Integrals, Series  and Other Formulas. CRC Press, Boca Raton, 2008.

\noindent [13] Apelblat, A., Kravitsky, N. Integral representations of derivatives and  integrals with respect to the order of the Bessel functions $J(t), I(t)$, the Anger  function $\textbf{J}(t)$ and the integral Bessel function $Ji(t)$. IMA J. Appl. Math\textit{.,}\textbf{ }1985\textbf{,  }34:187-210.

\noindent [14] Erd\'{e}lyi, A., Magnus, W., Oberhettinger, F., Tricomi, F.G. Tables of  Integral Transforms. McGraw-Hill, New York, 1954.

\noindent [15] Apelblat, A. Laplace Transforms and Their Applications. Nova Science  Publishers, Inc., New York, 2012.

\noindent [16] Apelblat, A. Volterra Functions. Nova Science Publishers, Inc., New York, 2008.

\noindent [17] Hladik, J. La Transformation de Laplace a Plusieurs Variables. Masson et Cie  \'{E}diteurs. Paris, 1969.

\noindent [18] Apelblat, A. The asymptotic limit of infinite integral of the Bessel function  $J_{\nu}(\nu x)$ as integral representation of elementary and special functions. Int. J. Appl.  Math. 2000, 2, 743-762.

\noindent [19] Apelblat, A. The application of the Dirac delta function $\delta(x - 1)$ to the evaluation of limits and integrals of elementary and special functions.
 Int. J. Appl. Math. 2004, 16, 323-339.

\noindent [20] Apelblat, A., The evaluation of the asymptotic relations, limits of series, integrals  and integral representations of elementary and special functions using the shifted  Dirac delta function $\delta(x - 1)$. Computing Letters (CoLe) 2008, 4, 11-19.

\noindent [21] Gradstein, I., Ryzhik, I. Tables of Series, Products and Integrals. Verlag Harri  Deutsch, Thun-Frankfurt, 1981.

\noindent [22] Abramowitz, M., Stegun, I.A. Handbook of Mathematical Functions with  Formulas, Graphs, and Mathematical Tables. U.S. National Bureau of Standards.  Applied Mathematics Series, 
 Washington, D.C. 1964.

\noindent[23] Garra, R., Giraldi, F., Mainardi,F.
Wright-type generalized coherent states,
WSEAS`Transactions on Mathematics (2019),  18, Art. 52, 428--431.

\noindent [24] Garra, R., Mainardi,F.
Some aspects of Wright functions in fractional differential equations,
Reports on Mathematical Physics (Poland) 2021,  87 (2),  265--273.
DOI: 10.1016/S0034-4877(21)00029-X,  E-print arXiv:2007.13340

\noindent [25] Mainardi, F., Consiglio, A.
 The Wright functions of the second kind in Mathematical Physics,
Mathematics (MDPI) 2020, Vol 8 N. 6, Art 884/1-26. 
DOI: 10.3390/MATH8060884; 
E-print arXiv:2007.02098 [math.GM] 


\noindent [26]  Li, C., Dao, X.
Fractional derivatives in complex planes. Nonlinear Anal. 71(5-6), 1857-1869, 2009.

\noindent [27] Guariglia, E.
 Fractional calculus, zeta functions and Shannon entropy, Open Mathematics, 19(1), 87-100, 2021.

\noindent [28] 
Ortigueira, M.D., Rodríguez-Germa, L., Trujillo, J.J.
 Complex Gr{\''u}nwald-Letnikov, Liouville, Riemann-Liouville, and Caputo derivatives for analytic functions, Commun. Nonlinear Sci. Numer. Simul. 16(11), 4174-4182, 2011.

\noindent [29] Guariglia, E.
Riemann zeta fractional derivative -- functional equation and link with primes, Advances in Difference Equations, 2019(1), 261/1-15, 2019.

\noindent [30] Z{\'a}vada, P.
 Operator of fractional derivative in the complex plane, 
 Comm. Math. Phys. 192, no. 2, 261–285, 1998.

\noindent [31] Lin, S-D, Srivastava, H.M. 
Some families of the Hurwitz–Lerch Zeta functions and associated fractional derivative and other integral representations. 
Appl. Math. Comput. 154(3), 725–733, 2004.

\noindent [32] Podlubny, I.
 Geometric  and physical interpretation of fractional integration and fractional differentiation. Fract. Calc. Appl. Anal. 5(4), 367–386,  2002.

\end{document}